\documentclass[11pt]{article}
\usepackage[a4paper,margin=28mm]{geometry}
\usepackage[T1]{fontenc}
\usepackage{lmodern,microtype}
\usepackage{amsmath,amssymb,amsthm,mathtools,mathrsfs}
\usepackage{enumitem,booktabs}
\usepackage[dvipsnames]{xcolor}
\usepackage[colorlinks=true,linkcolor=blue,citecolor=blue,urlcolor=blue]{hyperref}

\usepackage{fancyhdr}

\fancypagestyle{plain}{%
  \fancyhf{}%
  \fancyfoot[C]{\small\thepage}%
}

\numberwithin{equation}{section}
\newtheorem{theorem}{Theorem}[section]
\newtheorem{proposition}[theorem]{Proposition}
\newtheorem{lemma}[theorem]{Lemma}

\theoremstyle{definition}
\theoremstyle{remark}\newtheorem{remark}[theorem]{Remark}
\newcommand{\T}{\mathbb T}\newcommand{\Z}{\mathbb Z}
\newcommand{\R}{\mathbb R}\newcommand{\C}{\mathbb C}
\newcommand{\mes}{\operatorname{mes}}
\newcommand{\supp}{\operatorname{supp}}
\newcommand{\dd}{\,d}
\newcommand{\e}{\mathrm e}
\newcommand{\av}[1]{\langle #1\rangle}
\newcommand{\norm}[1]{\left\|#1\right\|}
\newcommand{\abs}[1]{\left|#1\right|}
\newcommand{\DC}{\mathrm{DC}}
\newcommand{\BMO}{\mathrm{BMO}}
\newcommand{\HH}{\mathcal H}

\newcommand{\eps}{\varepsilon}
\newcommand{\avg}[1]{\left\langle #1\right\rangle}

\newcommand{\one}{\mathbf 1}
\DeclareMathOperator{\dist}{dist}
\newcommand{\Tstrip}[1]{\mathbb T_{#1}}
\newcommand{\CJpi}{\mathcal C_{J,\pi}}
\newcommand{\AJpi}{\mathcal A_{J,\pi}}
\newcommand{\F}{\mathcal F}
\newcommand{\distT}[1]{\left\|#1\right\|_{\T}}
\newcommand{\1}{\mathbf 1}
\newcommand{\epss}{\varepsilon_*}
\newcommand{\gs}{\gamma_*}
\newcommand{\floor}[1]{\left\lfloor #1\right\rfloor}
\newcommand{\Csub}{C_{\mathrm{sub}}}
\allowdisplaybreaks[2]
\setlist[enumerate]{itemsep=3pt,topsep=5pt}

\title{Positivity and Continuity of Lyapunov Exponents for Multi-Frequency Skew-Shift Schr\"{o}dinger Operators }
\author{
Chao Wang\thanks{Email: \texttt{cwangmath@sohu.com}}
\qquad
Yuanyuan Peng\thanks{Email: \texttt{lunarpeng@foxmail.com}}
\qquad
Daxiong Piao\thanks{Corresponding author. Email: \texttt{dxpiao@ouc.edu.cn}}
\\[1.2ex]
\small School of Mathematical Sciences, Ocean University of China,\
\small Qingdao 266100, P.R. China
}

\date{}

\begin{document}
\maketitle

%--------------------------------------------------------------------------------------------
\begin{abstract}
We study discrete Schr\"odinger operators on \(\ell^2(\mathbb Z)\) generated by the vector-frequency two-step skew shift \(T_\omega(x,y)=(x+\omega,x+y)\), with an arbitrary nonconstant real-analytic sampling function on \(\mathbb T^d\times\mathbb T^d\). For frequencies in a fixed vector Diophantine class and sufficiently large coupling, we prove stretched-exponential large-deviation estimates at every sufficiently large integer scale, uniformly in the real energy and in the frequency. We also obtain positivity of the Lyapunov exponent for all real energies and a quantitative local continuity estimate in the energy. The proof addresses the simultaneous losses caused by the \(n^{-1}\) complex width induced by the shear, vector Diophantine small divisors in quadratic Weyl sums, and several-variable subharmonic concentration. The main ingredients are a vector Weyl estimate with Poisson regularization, a short-block large-deviation estimate on anisotropic shrinking strips, a several-variable boosting argument, and an all-offset avalanche comparison followed by strong induction.


\end{abstract}

\tableofcontents

% Section 1
\section{Introduction and main results}
\label{sec:intro}

\subsection{The model and the main theorem}

We consider a one-dimensional discrete Schr\"odinger operator whose
potential is sampled along a two-step affine transformation of a torus:
\begin{align}
 (H_{\lambda,\omega,z}\psi)(j)
 &=\psi(j+1)+\psi(j-1)+\lambda V(T_\omega^jz)\psi(j),
 \quad \psi\in\ell^2(\Z),\label{eq:originaloperator}\\
 T_\omega(x,y)&=(x+\omega,x+y),\quad
 z=(x,y)\in\T^d\times\T^d.
 \label{eq:originalmap}
\end{align}
Here $d\ge1$, $\lambda>0$, and $V:\T^{2d}\to\R$ is nonconstant and real
analytic. We refer to \eqref{eq:originalmap} as a
\emph{vector-frequency two-step skew shift}. The terminology is meant to
distinguish this model both from a multi-frequency torus rotation and
from higher-step skew shifts driven by a single scalar frequency: here
the nilpotent step remains two, while the frequency vector and both phase
variables are genuinely $d$-dimensional. Thus the base has dimension
$2d$, whereas the lattice remains one-dimensional and the transfer
matrices have size $2\times2$.

Throughout, torus measures are normalized Haar measures,
$e(t)=\exp(2\pi it)$, and $|k|=\sum_j|k_j|$. For $r>0$ write
$\T_r=\{z\in\C/\Z:|\operatorname{Im}z|<r\}$.
Choose $0<\rho<1/10$ so that $V$ has a bounded holomorphic extension to
$\T_\rho^{2d}$. The arithmetic assumption is the vector Diophantine
condition
\begin{equation}\label{eq:dc}
 \DC_d(\gamma,\tau)=\left\{\omega\in\T^d:
 \|k\cdot\omega\|_{\T}\ge\gamma |k|^{-\tau}
 \text{ for all }0\ne k\in\Z^d\right\},
 \qquad \gamma>0,\quad\tau>d.
\end{equation}
Define
\begin{align}
 A(z;\lambda,E)&=
 \begin{pmatrix}E-\lambda V(z)&-1\\1&0\end{pmatrix},
 \quad
 M_n(z;\lambda,E)
 =A(T_\omega^{n-1}z;\lambda,E)\cdots A(z;\lambda,E),
 \label{eq:originalmatrices}\\
 u_n(z;\lambda,E)&=
 n^{-1}\log\|M_n(z;\lambda,E)\|,
 \quad
 L_n(\lambda,E)=\int_{\T^{2d}}u_n(z;\lambda,E)\,dz.
 \label{eq:originalaverages}
\end{align}
The norm is the Euclidean operator norm. Subadditivity and invariance of
Haar measure give
\begin{equation}\label{eq:limit}
 L(\lambda,E)=\lim_{n\to\infty}L_n(\lambda,E)
 =\inf_{n\ge1}L_n(\lambda,E),\qquad L_{2n}\le L_n.
\end{equation}
Indeed $a_n=nL_n$ satisfies $a_{n+m}\le a_n+a_m$; for $N=qn+r$,
$0\le r<n$, division of $a_N\le qa_n+a_r$ by $N$ proves the limit
formula. Ergodicity, verified below, and Kingman's theorem \cite{Kingman}
also identify this number with almost-sure asymptotic matrix growth for
each fixed energy. We suppress the frequency argument in $L_n,L$.

Choose a fixed constant $C_V>1$, depending on the holomorphic bound for
$V$, large enough that
\begin{equation}\label{eq:S}
 S(\lambda,E)=\log\big(C_V(2+\lambda+|E|)\big)\ge1
\end{equation}
normalizes all one-step matrix bounds used below.

\begin{theorem}\label{thm:main}
Fix $d,\rho,\gamma,\tau$ as above and a nonconstant real-analytic $V$.
There exist $\lambda_0>1$, $N_0\ge1$, $C,c>0$, and
$b,\sigma\in(0,1)$ such that, for every
$\omega\in\DC_d(\gamma,\tau)$, $\lambda\ge\lambda_0$, $E\in\R$, and
integer $N\ge N_0$,
\begin{equation}\label{eq:mainLDT}
 \mes\left\{z:|u_N(z;\lambda,E)-L_N(\lambda,E)|
 >CS(\lambda,E)N^{-b}\right\}
 \le C\exp(-cN^\sigma).
\end{equation}
For the same parameters,
\begin{equation}\label{eq:mainLE}
 L(\lambda,E)\ge\tfrac12\log\lambda,\qquad E\in\R.
\end{equation}
For each fixed $\lambda\ge\lambda_0$ and compact interval $I\subset\R$,
there are $C_I,c_I>0$ such that
\begin{equation}\label{eq:continuitymain}
 |L(\lambda,E)-L(\lambda,E')|
 \le C_I\exp\big(-c_I|\log|E-E'||^\sigma\big),
 \quad E,E'\in I,\quad 0<|E-E'|<1/2.
\end{equation}
The constants in \eqref{eq:mainLDT}--\eqref{eq:mainLE} depend only on
$d,\rho,V,\gamma,\tau$, including the analytic norm and sublevel data of
$V$. The continuity constants may also depend on $\lambda,I$, and are
uniform in the same frequency class. The thresholds $\lambda_0,N_0$
are chosen independently of $E,N$ and of the particular frequency in
that class.
\end{theorem}

The modulus in \eqref{eq:continuitymain} is stretched-exponential in the
logarithm of the energy separation. It is stronger than every fixed
inverse power of that logarithm, but is not a power-law H\"older modulus.
The exceptional set in \eqref{eq:mainLDT} is defined for each energy
separately, with uniform constants; the theorem does not estimate its
union over all energies. The exponents are positive and explicit in the
proof, but are not optimized.

\begin{remark}[Arithmetic and frequency measure]
A scalar logarithmic Diophantine condition cannot be used in place of
\eqref{eq:dc} when $d\ge2$. The pigeonhole principle applied to the
$(Q+1)^d$ values $k\cdot\omega$ with $k\in\{0,\ldots,Q\}^d$ gives
$0<|k|_\infty\le Q$ with
$\|k\cdot\omega\|_{\T}\le C_dQ^{-d}$.
This contradicts a fixed lower bound
$\gamma/[|k|(1+\log|k|)^2]$ as $Q\to\infty$.
For the vector class, by contrast,
\[
 \mes(\T^d\setminus\DC_d(\gamma,\tau))
 \le2\gamma\sum_{k\ne0}|k|^{-\tau}
 \le C_{d,\tau}\gamma.
\]
Each nonzero integer character pushes Haar measure to Haar measure on
$\T$, and the lattice sum converges for $\tau>d$. In particular,
\[
 0<\gamma<(2C_{d,\tau})^{-1}
 \quad\text{gives}\quad
 \mes\DC_d(\gamma,\tau)\ge\tfrac12,
 \qquad
 \mes\left(\bigcup_{j\ge1}\DC_d(1/j,\tau)\right)=1.
\]
Theorem~\ref{thm:main} has a uniform coupling threshold on each fixed
class. On the full-measure union, the threshold may depend on the
frequency.
\end{remark}

For each prescribed $P'>0$, \eqref{eq:mainLDT} also gives
\[
 \mes\{|u_N-L_N|>CSN^{-b}\}\le N^{-P'}
 \qquad(N\ge N_{P'}),
\]
because $cN^\sigma-P'\log N-\log C\to+\infty$.
Only $N_{P'}$ depends on $P'$; the coupling threshold remains
$\lambda_0$.

\subsection{Related results and the role of the base dynamics}

Ergodic Schr\"odinger operators connect the growth of transfer matrices,
the spectral properties of a self-adjoint operator, and the recurrence
properties of the underlying dynamical system. General accounts are
given by Damanik \cite{Damanik} and Damanik--Fillman \cite{DF}.
Positivity of the Lyapunov exponent is a principal input in many
localization arguments, while quantitative deviations from finite-scale
means provide the control needed to compare different length scales.
For this reason, a uniform LDT is a natural objective even when
localization is not part of the theorem under consideration.

For analytic quasiperiodic sampling over a circle rotation, the
large-coupling positivity theorem of Sorets--Spencer \cite{SS} is a
basic reference. Bourgain--Goldstein \cite{BG} developed a
nonperturbative localization method based on large deviations and
elimination of resonances, while Goldstein--Schlag \cite{GS} introduced
the avalanche principle into quantitative regularity arguments for
Schr\"odinger cocycles. Bourgain--Jitomirskaya \cite{BJ} proved joint
continuity in frequency and energy in the one-frequency analytic
setting. In the higher-dimensional rotation setting, Bourgain
\cite{BourgainTd} established positivity and continuity of the
Lyapunov exponent for analytic sampling with arbitrary frequency
vectors. Goldstein--Schlag--Voda \cite{GSV} studied the large-coupling
spectrum of multi-frequency quasiperiodic Schr\"odinger operators, and
more recent continuity results for analytic multi-frequency cocycles
include Powell \cite{Powell} and Han--Zhang \cite{HanZhang}. These works
concern rotation bases; the statistical and complex-domain estimates
associated with a translation do not automatically survive the growing
shear in \eqref{eq:originalmap}. Bourgain's monograph
\cite{BourgainBook} gives a broader account of the multiscale framework
for deterministic lattice operators.

For the classical two-dimensional skew shift,
Bourgain--Goldstein--Schlag \cite{BGS} proved a large-coupling theorem
using quantitative averaging of subharmonic functions and the avalanche
principle; the localization argument additionally required recurrence
and resonance analysis. The preprint of Chan--Goldstein--Schlag
\cite{CGS} treats shifts and skew shifts with $C^\alpha$ potentials
under hypotheses adapted to low regularity. Klein \cite{Klein} proves
large-coupling positivity, continuity and localization for Gevrey
sampling under a transversality condition. In that work the skew-shift
model is two-dimensional, whereas its multi-frequency case is a
higher-dimensional torus translation. More recently, Liang--Fu
\cite{LiangFu} obtained a large-coupling asymptotic formula for the
Lyapunov exponent of Gevrey skew-shift Schr\"odinger operators. Tao
\cite{Tao2026} studied positivity and continuity of Lyapunov exponents
for discrete analytic Schr\"odinger equations with multivariable Gevrey
perturbations, including shift and skew-shift settings. These works are
close to the present problem in their large-coupling and regularity
aspects, but the vector-frequency two-step base
\eqref{eq:originalmap} and its associated anisotropic shrinking-strip
problem are different.

Quantitative work on the skew shift also includes the finite-scale
positivity criterion of Han--Lemm--Schlag \cite{HLSeffective} for the
cosine model with golden-ratio frequency and their analysis of quadratic
Weyl sums in a perturbative study of the Lyapunov exponent
\cite{HLSWeyl}. Kielstra--Lemm \cite{KielstraLemm} further studied the
finite-size Lyapunov exponent for the skew-shift potential, providing a
finite-scale perspective closely related to this program. Tao
\cite{Tao} proves nonperturbative positivity and related spectral
conclusions for skew-shift models under the hypotheses of that work.
Chen--Tao \cite{ChenTao} establish nonperturbative positivity and weak
H\"older continuity for a class of discrete multivariable Jacobi
operators over shifts and skew shifts. Our setting differs in that the
sampling function $V:\T^{2d}\to\R$ is an arbitrary fixed nonconstant
analytic function and the frequency belongs to a fixed vector
Diophantine class.

There is a second notion of ``higher-dimensional skew shift'' in which
the dimension is increased by increasing the nilpotent step. Kr\"uger
\cite{KrugerMSA} develops a multiscale approach to positivity for
ergodic operators, and Shi--Yuan \cite{ShiYuan} prove localization for a
class of Jacobi matrices associated with high-dimensional skew shifts.
This mechanism should be distinguished from
\eqref{eq:originalmap}: increasing the step changes the degree of the
polynomial orbit, whereas increasing $d$ here leaves the orbit
quadratic and instead introduces vector small divisors and
$2d$-dimensional Fourier summation. A separate spectral question is the
existence of intervals in the spectrum; Kr\"uger
\cite{KrugerIntervals} proves such a result in a semiclassical
skew-shift setting using additional spectral arguments, while
Goldstein--Schlag--Voda \cite{GSV} obtain spectral-interval results in a
large-coupling multi-frequency quasiperiodic setting. We establish the
large-deviation and Lyapunov-exponent conclusions for
\eqref{eq:originalmap}. The corresponding questions of Anderson
localization and the existence of intervals in the spectrum are pursued
separately.

The structural point relevant to the present paper is therefore the
following. The model \eqref{eq:originalmap} combines two features that
are usually separated in the references above: a genuinely vector
frequency and a genuine two-step skew shear. The step remains fixed, so
the orbit is still quadratic in time, but the Fourier analysis takes
place on $\Z^d\times\Z^d$ and must control the vector divisors
$k\cdot\omega$ and $\ell\cdot\omega$ uniformly on a fixed Diophantine
class. This distinction is the source of the additional quantitative
bookkeeping described next.

\subsection{Quantitative difficulties of the vector-frequency shear}

For the map \eqref{eq:originalmap}, direct induction gives
\begin{equation}\label{eq:introiterate}
 T_\omega^j(x,y)
 =
 \bigl(x+j\omega,\,
 y+jx+\tfrac12j(j-1)\omega\bigr),
 \qquad j\ge0.
\end{equation}
Thus the imaginary parts satisfy
\begin{equation}\label{eq:imaginaryshear}
 (\operatorname{Im}x,\operatorname{Im}y)
 =(\xi,\eta)
 \quad\longmapsto\quad
 (\xi,\eta+j\xi).
\end{equation}
The natural complex domain for the normalized logarithmic matrix norm
\[
 u_n(z;\lambda,E)
 :=
 \frac1n\log\|M_n(z;\lambda,E)\|
\]
therefore has width of order $n^{-1}$ in the $x$-variables and order
one in the $y$-variables. More precisely, the shear in
\eqref{eq:imaginaryshear} forces an anisotropic complex neighborhood of
the form
\[
 \|\operatorname{Im}x\|_\infty\lesssim n^{-1},
 \qquad
 \|\operatorname{Im}y\|_\infty\lesssim1.
\]
Moreover, $u_n$ is plurisubharmonic rather than holomorphic.
Consequently, analyticity of the sampling function alone does not
produce a uniform exponential Fourier tail for $u_n$.

The multi-frequency aspect adds a separate arithmetic difficulty.
For a Fourier character indexed by
$(k,\ell)\in\Z^d\times\Z^d$, the orbit phase is
\[
 (k,\ell)\cdot T_\omega^j(x,y)
 =
 k\cdot x+\ell\cdot y
 +j(\ell\cdot x+k\cdot\omega)
 +\frac{j(j-1)}2\,\ell\cdot\omega.
\]
A van der Corput difference therefore produces linear frequencies
$h\ell\cdot\omega$. When $d>1$, these are genuinely vector small
divisors, and the estimate must be uniform for all
$\ell\in\Z^d\setminus\{0\}$ and
$\omega\in\DC_d(\gamma,\tau)$.

The arithmetic and complex-analytic losses also interact quantitatively.
After Poisson regularization, the vector Weyl estimate is summed over
$\mathbf m=(k,\ell)\in\Z^{2d}$ and produces the bound
\[
 R^{-\delta}
 \sum_{\mathbf m\ne0}
 e^{-t|\mathbf m|_1}(1+|\mathbf m|_1)^\tau
 \lesssim
 R^{-\delta}t^{-(2d+\tau)}.
\]
On the other hand, the shrinking $x$-strip gives a one-sided
subharmonic regularization cost of order $n^4t$. Thus the short-block
argument has to make the two losses
\[
 n^4t
 \qquad\text{and}\qquad
 R^{-\delta}t^{-(2d+\tau)}
\]
simultaneously small. This is the quantitative reason for imposing a
polynomial relation $R\ge R_0n^A$ between the orbit-averaging scale and
the block length.

There is one further complex-domain issue. The full orbit average
\[
 W_R(z)=\frac1R\sum_{j=0}^{R-1}w_n(T_\omega^jz)
\]
is not used through a uniformly bounded extension to the original
product strip; such an extension is not available at the required
scale because the shear amplifies the imaginary $x$-component.
Instead, the proof uses a sequential complexification. It first averages
in the sheared $y$-variables to obtain a legitimate complex extension
of a function of $x$, and then, for real $x$, treats the corresponding
$y$-sections. This triangular use of the complex variables is what
allows one-sided subharmonic control to be converted into a two-sided
short-block deviation estimate.

These features separate the present problem into four quantitative
tasks:
\[
\boxed{
\begin{aligned}
&\text{shrinking-strip subharmonic concentration},\\
&\text{vector Diophantine cancellation for quadratic orbit phases},\\
&\text{short-block averaging compatible with the shear},\\
&\text{propagation from finite scales to every large integer scale}.
\end{aligned}}
\]

\subsection{Proof architecture and contributions}

The proof closes these tasks through the following quantitative
ingredients. Each is established below in the form needed for the
argument, rather than assumed through an operator-specific
large-deviation theorem.

\begin{enumerate}[label=(\arabic*)]
\item \emph{Several-variable boosting on shrinking strips.}
For a continuous separately subharmonic function $f$ on a product strip
of width $\rho/N$, Proposition~\ref{prop:boosting} shows that a real
decomposition
\[
 f=f_0+f_1,\qquad
 \|f_0\|_\infty\le N^{-s},\qquad
 \|f_1\|_1\le N^{-P_m(s)}
\]
implies
\[
 \mes\{|f|>CN^{-\beta_m(s)}\}
 \le Ce^{-cN^{\beta_m(s)}}.
\]
The exponent $P_m(s)$ is chosen recursively. The proof isolates the
strip loss in a one-dimensional BMO estimate and then combines
coordinatewise Fej\'er smoothing with a finite net. The classical
harmonic-analysis inputs are those of \cite{JN,Stein,Zygmund}; what is
retained here is explicit control of the loss under dimension induction
and shrinking strip width.

\item \emph{A vector-frequency short-block LDT.}
For $|w_n|\le1$ on
$\T^d_{\rho_x/n}\times\T^d_{\rho_y}$,
Proposition~\ref{prop:short-block-LDT} gives, whenever
$R\ge R_0n^A$,
\[
 \mes\left\{
 \left|R^{-1}\sum_{j<R}w_n\circ T^j-\langle w_n\rangle\right|
 >CR^{-a}
 \right\}
 \le Ce^{-cR^a}.
\]
The proof combines the vector Weyl bound with Poisson regularization and
then uses the sequential complex-domain argument described above. In
particular, it does not require a uniformly bounded extension of the
full orbit average to the original product strip.

\item \emph{A uniform large-coupling initial scale for arbitrary analytic
sampling.}
The normalization of $u_n$ on the anisotropic strip is combined with a
uniform analytic sublevel estimate
\[
 \mes\{|V-t|<\varepsilon\}
 \le C_{\rm sub}\varepsilon^\beta,
 \qquad t\in\R.
\]
No product structure or coordinatewise nondegeneracy is imposed on
$V$. Outside the resulting small bad set, a finite-volume determinant
factorization gives a pointwise lower bound for the transfer matrix; for
energies outside the coupling scale, direct diagonal dominance gives the
same conclusions without an exceptional set. This supplies the initial
LDT, positivity and the first two-scale gap estimate uniformly in the
real energy.

\item \emph{Propagation to every sufficiently large integer scale.}
Averaging the avalanche principle over all block offsets gives an error
of order $n/N$. The resulting bad-set bound $CNq_n$ enters the boosting
estimate as a small $L^1$ remainder. With the parent-scale rule
\[
 n=\lfloor N^{1/A}/K\rfloor,
\]
a strong induction over the integers propagates both the LDT and the
positive lower bound for the Lyapunov exponent. The coupling threshold
is fixed once at the initial scale, and the cumulative loss in the lower
bound is summable.
\end{enumerate}

The new point is the quantitative closure of these mechanisms for a
genuinely vector-frequency two-step skew shift. Three losses must be
controlled simultaneously: the $n^{-1}$ complex width created by the
shear, the vector Diophantine losses in the quadratic Weyl sums, and the
dimension-dependent losses in several-variable subharmonic
concentration. The proof keeps these losses compatible through the
short-block averaging scale and then propagates the resulting estimates
to every sufficiently large integer scale.

The individual ingredients themselves have classical precedents: the
avalanche principle, BMO and John--Nirenberg estimates, Poisson
regularization, and analytic sublevel bounds are established tools. The
contribution lies in combining them in a form that preserves the
dependence on strip width, dimension, scale and the fixed vector
Diophantine class. The sampling function may couple all $2d$ phase
variables or factor through fewer coordinates; no artificial
coordinatewise interaction is assumed. The resulting stretched-
exponential exceptional-set estimate also implies every prescribed
polynomial exceptional-set bound with the same fixed large-coupling
threshold. Its use in localization requires additional energy and
recurrence quantifiers and is not part of Theorem~\ref{thm:main}.

\subsection{Notation, constants, and organization}

We use $H,A,M_n,u_n$ as defined in
\eqref{eq:originaloperator}--\eqref{eq:originalaverages} throughout,
and write $T=T_\omega$ when the frequency is fixed. The transfer
convention is
\begin{equation}\label{eq:transferconvention}
 \binom{\psi(n)}{\psi(n-1)}
 =M_n(z;\lambda,E)
 \binom{\psi(0)}{\psi(-1)},
 \qquad H\psi=E\psi.
\end{equation}
All averages and exceptional-set measures below refer to the real torus.

\paragraph{Constant convention.}
In the one-dimensional estimates, numerical constants are absolute
apart from explicitly specified bounds. In
Section~\ref{sec:subharmonic}, the constants and lower scale thresholds
depend only on the stated dimension, strip parameters, decomposition
exponent and absolute bounds. In
Section~\ref{sec:skew-averaging} they may also depend on
$\gamma,\tau$. In Sections~\ref{sec:initial}--\ref{sec:main-proof} the
fixed data are $(d,\rho,V,\gamma,\tau)$; once these and the auxiliary
exponents are fixed, all lower thresholds are chosen uniformly in
$\omega\in\DC_d(\gamma,\tau)$, $E\in\R$ and
$\lambda\ge\lambda_0$. Increasing a fixed threshold never means changing
it with the current scale or with the input function. The order of the
scale and coupling choices is specified in
Section~\ref{sec:main-proof}.

\paragraph{Organization.}
Section~\ref{sec:subharmonic} develops the one-dimensional strip
estimates. Section~\ref{sec:boosting} proves several-variable boosting
and one-sided concentration on shrinking strips.
Section~\ref{sec:skew-averaging} establishes the vector Weyl estimate
and the short-block LDT. Section~\ref{sec:initial} constructs the
large-coupling initial scale, and Section~\ref{sec:AP} gives the
all-offset avalanche comparison and the one-step propagation estimate.
Finally, Section~\ref{sec:main-proof} closes the strong induction over
all sufficiently large integer scales and proves positivity and
continuity in the energy.

% Section 2
\section{One-dimensional subharmonic estimates with strip dependence}
\label{sec:subharmonic}

We prove the precise analytic input, including its dependence on strip
width. The harmonic-analysis facts used are the periodic Hilbert transform
map $L^\infty\to\BMO$ and the John--Nirenberg inequality; see \cite{JN,Stein}.
They have absolute constants under our normalization of $\T$.

\subsection{Notation and harmonic-analysis preliminaries}
%\paragraph{The torus, integration, and Fourier coefficients.}
We identify $\T=\R/\Z$ with $[0,1)$ and use normalized Lebesgue
measure $dx$. The distance to the nearest integer is
$\|x\|_\T=\min_{j\in\Z}|x-j|$. For $f\in L^1(\T)$, set
\begin{equation}\label{eq:Fourierdefinitions}
 e(t)=e^{2\pi it},\qquad
 \widehat f(k)=\int_0^1f(x)e(-kx)\,dx,\qquad
 \langle f\rangle=\widehat f(0),\quad k\in\Z.
\end{equation}
The symbols $\|f\|_p$ and $\|f\|_\infty$ denote the usual $L^p$ norm
and essential supremum on the indicated real torus; for a continuous
function the latter is its supremum. For $g\in L^1(\T)$, convolution is
\begin{equation}\label{eq:convolutiondefinition}
 (g*f)(x)=\int_0^1g(x-t)f(t)\,dt,
 \qquad \widehat{g*f}(k)=\widehat g(k)\widehat f(k).
\end{equation}
For a finite measure $\nu$ the same notation means
$(\nu*f)(x)=\int_\T f(x-t)\,d\nu(t)$. We write $\delta_0$ for the unit
point mass at zero, $\mathbf1_E$ for the indicator of a set $E$, and
$a_+=\max(a,0)$. On $\T^m$, the same conventions apply with
$k\cdot x$ in place of $kx$ and $|k|_1=\sum_{j=1}^m|k_j|$.

%\paragraph{Subharmonicity and distributions.}
For $z=x+iy$, let $dA(z)=dx\,dy$,
$
D(z,r)=\{w:|w-z|<r\},
$
and $\Delta=\partial_x^2+\partial_y^2$. A continuous real-valued
function $f$ is subharmonic if, for every closed disk in its domain,
\begin{equation}\label{eq:submeanconvention}
 f(z)\le\frac1{2\pi}\int_0^{2\pi}f(z+re^{i\theta})\,d\theta.
\end{equation}
It is separately subharmonic on a product domain if this condition holds
in each complex coordinate with the others fixed. For a locally integrable
function, distributional differentiation means
\[
 \langle\Delta f,\varphi\rangle=\int f\Delta\varphi\,dA,
 \qquad \varphi\in C_c^\infty,
\]
where $C_c^\infty$ is the space of smooth compactly supported test
functions. The Riesz measure and kernel normalization are
\begin{equation}\label{eq:Rieszconventions}
 \mu_f=\frac1{2\pi}\Delta f\ge0,
 \qquad \Delta\log|z|=2\pi\delta_0.
\end{equation}
These conventions apply on local lifts of the cylinder $\C/\Z$;
see \cite{Ransford,Hormander}. An inequality between distributions is
tested against nonnegative $\varphi$.

%\paragraph{The Poisson kernel and the operator $P_t$.}
The radial parameter $a\in[0,1)$ and the semigroup parameter $t>0$
are related by $a=e^{-t}$. Define
\begin{align}
 \mathscr P_a(x)&=\frac{1-a^2}{1-2a\cos(2\pi x)+a^2}
     =\sum_{k\in\Z}a^{|k|}e(kx),\label{eq:Poissonkernel}\\
 P_tf(x)&=(\mathscr P_{e^{-t}}*f)(x)
     =\int_0^1\mathscr P_{e^{-t}}(x-u)f(u)\,du.
     \label{eq:Poissonoperator}
\end{align}
Thus $P_t$ is an averaging operator, not a scalar or a pointwise
multiplier. The identities used below are
\begin{align}
 &\mathscr P_a\ge0,\qquad \int_\T\mathscr P_a=1,
       \qquad \widehat{P_tf}(k)=e^{-t|k|}\widehat f(k),
       \label{eq:Poissonproperties}\\
 &P_t1=1,\quad \langle P_tf\rangle=\langle f\rangle,
       \quad P_tP_s=P_{t+s},\quad
       \|P_tf\|_p\le\|f\|_p\quad(1\le p\le\infty).
       \label{eq:Poissoncontraction}
\end{align}
Moreover $P_tf\to f$ in $L^p$ for $1\le p<\infty$, and uniformly
for continuous $f$, as $t\downarrow0$. In the weak measure sense,
$\mathscr P_a(x)\,dx\to\delta_0$ as $a\uparrow1$; we write
$\mathscr P_1=\delta_0$. All these statements follow from the Fourier
series and the positive approximate-identity property
\cite{Zygmund,SteinShakarchi}.

In $m$ variables let $P_t^{(j)}$ act in coordinate $j$ alone and set
\begin{equation}\label{eq:productPoisson}
 Q_t=\prod_{j=1}^mP_t^{(j)},\qquad
 \widehat{Q_tf}(k)=e^{-t|k|_1}\widehat f(k).
\end{equation}
The factors commute. Positive real convolution preserves the complex
strip and separate subharmonicity, as proved in
Lemma~\ref{lem:positive-average}. For bounded $f$, the Fourier series
of $Q_tf$ converges absolutely at every $t>0$.

%\paragraph{Hilbert transform and mean oscillation.}
The periodic Hilbert transform $\HH$ is defined on $L^2(\T)$ by
\begin{equation}\label{eq:Hilbertdefinition}
 \widehat{\HH f}(k)=-i\operatorname{sgn}(k)\widehat f(k),
 \qquad \operatorname{sgn}(0)=0.
\end{equation}
It extends to $L^p$, $1<p<\infty$, and to distributions by the same
multiplier. On smooth functions it equals the principal-value convolution
with $\cot(\pi x)$. In particular,
\begin{equation}\label{eq:Hilbertidentities}
 \HH^2f=-f+\langle f\rangle,\quad
 \langle\HH f\rangle=0,\quad
 \HH P_t=P_t\HH.
\end{equation}
For a real integrable $f$, define its BMO seminorm by
\begin{equation}\label{eq:BMOdefinition}
 \|f\|_{\BMO(\T)}=\sup_{I\subset\T}\frac1{|I|}
       \int_I|f-f_I|\,dx,
 \qquad f_I=\frac1{|I|}\int_I f\,dx.
\end{equation}
The supremum runs over arcs of positive length, including the whole
circle. Constants have seminorm zero. The two estimates used below are
\begin{align}
 \|\HH g\|_{\BMO(\T)}&\le C\|g\|_\infty,
       \label{eq:HilbertBMO}\\
 \mes\{|f-\langle f\rangle|>u\}
   &\le C\exp\left(-\frac{cu}{\|f\|_{\BMO(\T)}}\right),\qquad u>0.
       \label{eq:JN}
\end{align}
These are the periodic $L^\infty$-to-BMO estimate and the
John--Nirenberg inequality \cite{Stein,JN}. In \eqref{eq:JN}, the
zero-seminorm case means that the exceptional set has measure zero.
Only one-dimensional versions of these estimates are used.

\subsection{Riesz mass and the harmonic remainder}
We use the distributional and representative conventions of classical
potential theory; see \cite{Ransford,Hormander}. The following proof
tracks their constants on a varying periodic strip.

\begin{lemma}\label{lem:riesz}
Suppose $0<r<1/10$, $f$ is continuous and subharmonic on $\T_r$, and $|f|\le 1$. Then on the real circle it has a representation
\begin{equation}\label{eq:riesz}
f(x)=h(x)+c_0+\int \log|1-a_\zeta e(x-b_\zeta)|\,d\nu(\zeta),
\qquad 0\le a_\zeta\le 1.
\end{equation}
where $\nu$ is a positive measure supported in $\T\times(-r/2,r/2)$,
\[
\nu\bigl(\T\times(-r/2,r/2)\bigr)\le \frac{C}{r},
\]
and $h$ is harmonic on $\T_{r/2}$. On $\T_{r/4}$ one has
\[
|h|\le \frac{C}{r^2}.
\]
Consequently, for $k\ne 0$,
\[
|\widehat h(k)|\le Cr^{-2}e^{-cr|k|},\qquad
\sum_{k\ne 0}|k|\,|\widehat h(k)|\le Cr^{-4},
\]
and
\[
|\widehat f(k)|\le Cr^{-4}|k|^{-1}.
\]
Here $e(t)=\exp(2\pi i t)$.
\end{lemma}

\begin{proof}
Since $f$ is subharmonic, by (\ref{eq:Rieszconventions}), we have $\mu:=\mu_f$ is positive. Fix $\chi_0\in C_c^\infty((-1,1))$ such that
\[
0\le \chi_0\le 1,
\qquad
\chi_0\equiv 1\quad\text{on }[-3/4,3/4],
\]
and set
\[
\chi_r(y):=\chi_0(y/r).
\]
Then
\begin{equation}\label{eq:cutoff}
\chi_r(y)=1\quad (|y|\le 3r/4),\qquad
\supp\chi_r\subset(-r,r),\qquad
|\chi_r''(y)|\le Cr^{-2}.
\end{equation}
Since $\chi_r$ depends only on $y=\operatorname{Im}z$, we have
$\Delta\chi_r=\chi_r''$. Positivity of $\mu$ and distributional
integration by parts give
\begin{align}
\mu\bigl(\T\times[-3r/4,3r/4]\bigr)
&\le \int_{\T_r}\chi_r\,\dd\mu \notag\\
&=\frac{1}{2\pi}\int_{-r}^{r}\int_0^1
f(x+iy)\chi_r''(y)\,\dd x\,\dd y.\label{eq:mass1}
\end{align}
Using $|f|\le 1$ and \eqref{eq:cutoff},
\begin{align}
\int_{\T_r}\chi_r\,\dd\mu
&\le \frac{1}{2\pi}\int_{-r}^{r}\int_0^1
|\chi_r''(y)|\,\dd x\,\dd y \notag\\
&\le Cr^{-2}(2r)
\le \frac{C}{r}.\label{eq:mass2}
\end{align}
Hence, with
$
\nu:=\mu\big|_{\T\times(-r/2,r/2)},
$
we have
\begin{equation}\label{eq:numass}
\nu\bigl(\T\times(-r/2,r/2)\bigr)\le \frac{C}{r}.
\end{equation}

Set
\begin{equation}\label{eq:potential}
K(z):=\log|2\sin \pi z|,
\qquad
U(z):=\int_{\T\times(-r/2,r/2)}K(z-\zeta)\,\dd\nu(\zeta).
\end{equation}
The kernel $K$ is locally integrable on $\C/\Z$. Indeed, near its
singularity it differs from $\log|z|$ by a bounded function, and
\[
\int_{|z|<\varepsilon}|\log|z||\,\dd A(z)
=2\pi\int_0^\varepsilon |\log \rho|\,\rho\,\dd\rho<\infty.
\]
Consequently, for every compact $E\Subset\T_{r/2}$,
\begin{align*}
\int_E |U(z)|\,\dd A(z)
&\le \int
\left(\int_E|K(z-\zeta)|\,\dd A(z)\right)\dd\nu(\zeta)\notag\\
&\le C_E\,\nu\bigl(\T\times(-r/2,r/2)\bigr)<\infty,
\end{align*}
so $U\in L^1_{\mathrm{loc}}(\T_{r/2})$.

On the periodic cylinder, $\Delta K=2\pi\delta_0$ in the sense of distributions. Hence, for
$\varphi\in C_c^\infty(\T_{r/2})$,
\[
\langle \Delta U,\varphi\rangle =\int
\left\langle \Delta K(\cdot-\zeta),\varphi\right\rangle
\dd\nu(\zeta)=2\pi\int \varphi(\zeta)\,\dd\nu(\zeta),
\]
and therefore
\begin{equation}\label{eq:lapU}
\Delta U=2\pi\nu.
\end{equation}
Since $\nu=\mu$ on the open strip $\T_{r/2}$,
\[
\Delta(f-U)=2\pi(\mu-\nu)=0
\qquad\text{on }\T_{r/2}.
\]
By Weyl's lemma, $f-U$ agrees almost everywhere with a smooth harmonic
function $h$ on $\T_{r/2}$; see, for example, \cite{Hormander}.
The functions $f$ and $h+U$, taken with their standard subharmonic
representatives, agree almost everywhere and hence everywhere. Equivalently,
for every disk $D(z,\rho)\Subset\T_{r/2}$,
\begin{align}
f(z)
&=\lim_{\rho\downarrow0}\frac{1}{\pi\rho^2}
\int_{D(z,\rho)}f(w)\,\dd A(w)\notag\\
&=\lim_{\rho\downarrow0}\frac{1}{\pi\rho^2}
\int_{D(z,\rho)}(h(w)+U(w))\,\dd A(w)\notag\\
&=h(z)+U(z).\label{eq:pointwise}
\end{align}
Thus
\begin{equation}\label{eq:fhU}
f=h+U\qquad\text{pointwise on }\T_{r/2}.
\end{equation}
In particular, $U=f-h$ is continuous there.

For $z=x+iy\in\T_{r/2}$ and
$\zeta=b+i\eta\in\supp\nu$, one has $|y-\eta|<r$. Hence
\[
|2\sin\pi(z-\zeta)|^2 =4\Bigl(\sin^2\pi(x-b)+\sinh^2\pi(y-\eta)\Bigr)\le 4\cosh^2(\pi r),
\]
so, since $r<1/10$,
\[
K(z-\zeta)\le C.
\]
By \eqref{eq:numass},
\[
U(z)\le C\nu\bigl(\T\times(-r/2,r/2)\bigr)
\le \frac{C}{r},
\qquad z\in\T_{r/2}.
\]
Since $f\ge -1$, \eqref{eq:fhU} gives
\begin{equation}\label{eq:hlower}
h(z)=f(z)-U(z)\ge -1-\frac{C}{r}\ge -\frac{C}{r},
\qquad z\in\T_{r/2}.
\end{equation}

On the real circle, if $\zeta=b+i\eta$, then
\begin{equation}\label{eq:kernel-factor}
K(x-\zeta)
=\pi|\eta|+\log|1-e^{-2\pi|\eta|}e(x-b)|.
\end{equation}
Set
\begin{equation}\label{eq:ab}
a_\zeta:=e^{-2\pi|\eta|},
\qquad
b_\zeta:=b\pmod 1,
\qquad
c_0:=\pi\int |\operatorname{Im}\zeta|\,\dd\nu(\zeta).
\end{equation}
Then $0<a_\zeta\le 1$, and \eqref{eq:fhU}, \eqref{eq:potential}, and
\eqref{eq:kernel-factor} give precisely \eqref{eq:riesz}. Moreover,
\[
0\le c_0\le \frac{\pi r}{2}\nu(\T_{r/2})\le C.
\]

For $|y|<r/2$ and $\zeta=b+i\eta$, put $t=y-\eta$. We claim that
\begin{equation}\label{eq:horizontalmean}
\int_0^1K(x+iy-\zeta)\,\dd x=\pi|t|.
\end{equation}
Indeed, for $t>0$,
\[
|2\sin\pi(x+it)|
=e^{\pi t}|1-e^{-2\pi t}e(x)|,
\]
whence
\[
K(x+it)=\pi t+\log|1-e^{-2\pi t}e(x)|.
\]
Since $0<e^{-2\pi t}<1$,
\[
\int_0^1\log|1-e^{-2\pi t}e(x)|\,\dd x=0,
\]
and therefore the left-hand side of \eqref{eq:horizontalmean} equals
$\pi t$. The case $t<0$ is identical, while $t=0$ follows by
$L^1$-continuity. Thus \eqref{eq:horizontalmean} holds.

By the uniform \(L^1\)-integrability of the logarithmic kernel, Fubini's theorem, \eqref{eq:horizontalmean}, and \eqref{eq:fhU},
\begin{equation}\label{eq:hmean}
\int_0^1 h(x+iy)\,\dd x=\int_0^1 f(x+iy)\,\dd x-\int \pi|y-\operatorname{Im}\zeta|\,\dd\nu(\zeta)\le 1,
\qquad |y|<r/2.
\end{equation}
Choose an absolute $C_1>0$ so large that, by \eqref{eq:hlower},
\[
H(z):=h(z)+\frac{C_1}{r}+2>0,
\qquad z\in\T_{r/2}.
\]
Moreover, \eqref{eq:hmean} implies
\begin{equation}\label{eq:Hmean}
\int_0^1H(x+iy)\,\dd x
\le 1+\frac{C_1}{r}+2
\le \frac{C}{r},
\qquad |y|<r/2.
\end{equation}
If $z\in\T_{r/4}$, then
$D(z,r/8)\Subset\T_{3r/8}\subset\T_{r/2}$. Since $H$ is positive and
harmonic, the mean-value formula and \eqref{eq:Hmean} yield
\[
H(z)
=\frac{64}{\pi r^2}\int_{D(z,r/8)}H(w)\,\dd A(w)\notag\\
\le \frac{64}{\pi r^2}\int_{-3r/8}^{3r/8}
\int_0^1 H(x+iy)\,\dd x\,\dd y\notag\\
\le \frac{C}{r^2}\int_{-3r/8}^{3r/8}\frac{1}{r}\,\dd y
\le \frac{C}{r^2}.
\]
Together with \eqref{eq:hlower} and $0<r<1$, this gives
\begin{equation}\label{eq:hstrip}
|h(z)|\le \frac{C}{r^2},
\qquad z\in\T_{r/4}.
\end{equation}

For $k\ne0$, define
\[
h_k(y):=\int_0^1 h(x+iy)e(-kx)\,\dd x.
\]
Since $\Delta h=0$, two integrations by parts in the periodic variable
$x$ give
\begin{align}
h_k''(y)
&=\int_0^1 h_{yy}(x+iy)e(-kx)\,\dd x\notag\\
&=-\int_0^1 h_{xx}(x+iy)e(-kx)\,\dd x\notag\\
&=(2\pi k)^2h_k(y).\label{eq:modeode}
\end{align}
Writing $a=r/8$, the general solution of \eqref{eq:modeode} gives
\[
h_k(0)=\frac{h_k(a)+h_k(-a)}{2\cosh(2\pi|k|a)}.
\]
By \eqref{eq:hstrip},
\[
|h_k(\pm r/8)|\le Cr^{-2},
\]
and therefore
\[
|\widehat h(k)|=|h_k(0)|
\le Cr^{-2}e^{-cr|k|}.
\]
Consequently,
\[
\sum_{k\ne0}|k|\,|\widehat h(k)|
\le Cr^{-2}\sum_{k\ne0}|k|e^{-cr|k|}
\le Cr^{-2}\frac{e^{-cr}}{(1-e^{-cr})^2}
\le Cr^{-4},
\]
where the last inequality uses $1-e^{-cr}\ge c_1r$ for
$0<r<1/10$.

It remains to estimate the Fourier coefficients of $U$. For
$0\le a<1$,
\[
\log|1-ae(x-b)|
=-\sum_{m=1}^{\infty}\frac{a^m}{m}
\cos(2\pi m(x-b)),
\]
so for $k\ne0$,
\[
\widehat{\log|1-ae(\cdot-b)|}(k)
=-\frac{a^{|k|}}{2|k|}e(-kb).
\]
For $a=1$, the same identity follows by letting $a\uparrow1$ in
$L^1(\T)$, indeed $\log|1-e(x-b)|\in L^1(\T)$. Hence, using
\eqref{eq:numass},
\[
|\widehat U(k)|
\le \int\frac{a_\zeta^{|k|}}{2|k|}\,\dd\nu(\zeta)
\le \frac{1}{2|k|}\nu(\T_{r/2})
\le \frac{C}{r|k|}.
\]
Since $f=h+U$ on the real circle,
\[
|\widehat f(k)|
\le Cr^{-2}e^{-cr|k|}+\frac{C}{r|k|}.
\]
Finally, with $t=r|k|$ and $te^{-ct}\le C$,
\[
r^{-2}e^{-cr|k|}
\le \frac{Cr^{-3}}{|k|}.
\]
Since $0<r<1$,
\[
|\widehat f(k)|
\le \frac{Cr^{-3}+Cr^{-1}}{|k|}
\le \frac{Cr^{-4}}{|k|}.
\]
This proves all assertions.
\end{proof}

\subsection{One-sided Poisson regularization and the BMO estimate}
\label{subsec:poisson-bmo}

Lemma \ref{lem:riesz} gives a polynomial loss in the strip width.  Throughout this
subsection we fix
\begin{equation}
    K_r=C_*r^{-4}\ge 1,
    \label{eq:pzero}
\end{equation}
where the absolute constant $C_*$ is chosen sufficiently large to dominate
all constants in the Riesz-mass and harmonic Fourier estimates of
Lemma \ref{lem:riesz}. % Let $P_t$ be the periodic Poisson semigroup,
%\begin{equation}
%    \widehat{P_tg}(k)=e^{-t|k|}\widehat g(k),
%    \qquad t>0.
%    \label{eq:Pt-def}
%\end{equation}
%No optimality of the exponent $4$ in \eqref{eq:pzero} will be used.

\begin{lemma}
\label{lem:oneD}
Assume the hypotheses of Lemma \ref{lem:riesz}. Then, for $0<t<1$,
\begin{equation}
    f(x)\le P_tf(x)+K_rt,
    \qquad x\in\T.
    \label{eq:one-sided-poisson}
\end{equation}
Suppose, on the real circle,
\begin{equation}
    f=f_0+f_1,
    \qquad
    \norm{f_0}_{L^\infty(\T)}\le \eps_0,
    \qquad
    \norm{f_1}_{L^1(\T)}\le \eps_1.
    \label{eq:bmo-decomp-assumption}
\end{equation}
Then, for every $L\ge2$,
\begin{equation}
    \norm{f}_{\BMO(\T)}
    \le
    C\left(
        \eps_0\log L
        +\eps_1L
        +\frac{K_r}{L}
    \right).
    \label{eq:bmo-decomp}
\end{equation}
In particular, if
\begin{equation}
    \norm{f-\avg f}_{L^1(\T)}\le\eps\le1,
    \label{eq:l1-smallness}
\end{equation}
then
\begin{equation}
    \norm{f}_{\BMO(\T)}
    \le C\sqrt{K_r\eps}.
    \label{eq:bmo-l1}
\end{equation}
The same conclusions hold under any fixed bound $|f|\le B$, with the
constant depending on $B$.
\end{lemma}

\begin{proof}
For $0\le a\le1$ and $b\in\T$, set
$
\ell_{a,b}(x):=\log|1-ae(x-b)|.
$
For $k\ne0$,
\begin{equation}
    \widehat{\ell_{a,b}}(k)
    =-
    \frac{a^{|k|}}{2|k|}e(-kb),
    \qquad
    \widehat{\ell_{a,b}}(0)=0,
    \label{eq:ellab-fourier}
\end{equation}
where the endpoint $a=1$ is understood in $L^1(\T)$.  Put $q=e^{-t}$, by \eqref{eq:Poissonproperties} and \eqref{eq:ellab-fourier},
\[
    P_t\ell_{a,b}=\ell_{qa,b}.
\]
For $\theta\in\R$,
\begin{align}
    |1-qae^{i\theta}|^2-q|1-ae^{i\theta}|^2
    &=1-2qa\cos\theta+q^2a^2
      -q(1-2a\cos\theta+a^2) \notag\\
    &=(1-q)(1-qa^2)\ge0.
    \label{eq:elementary-poisson-ineq}
\end{align}
Hence
$
    |1-qae^{i\theta}|^2
    \ge q|1-ae^{i\theta}|^2,
$
and therefore
\begin{equation}
    \ell_{a,b}(x)
    \le \ell_{qa,b}(x)+\frac t2
    =P_t\ell_{a,b}(x)+\frac t2.
    \label{eq:log-one-sided}
\end{equation}

Write the Riesz representation of Lemma~\ref{lem:riesz} as
\begin{equation}
    f=h+c_0+\int \ell_{a_\zeta,b_\zeta}\,d\nu(\zeta).
    \label{eq:Riesz-used}
\end{equation}
For the harmonic term,
$
    h-P_th=\sum_{k\ne0}
      \bigl(1-e^{-t|k|}\bigr)\widehat h(k)e(kx),
$
so, using $1-e^{-u}\le u$ for $u\ge0$,
\[
    \norm{h-P_th}_{L^\infty(\T)}
    \le
    \sum_{k\ne0}
    \bigl(1-e^{-t|k|}\bigr)|\widehat h(k)|
    \le
    t\sum_{k\ne0}|k|\,|\widehat h(k)|
    \le C r^{-4}t.
\]
Applying \eqref{eq:log-one-sided} in \eqref{eq:Riesz-used},
\[
    f-P_tf
    \le
    C r^{-4}t+\frac t2\nu(\T_{r/2})
    \le
    C(r^{-4}+r^{-1})t
    \le K_rt,
\]
which proves \eqref{eq:one-sided-poisson}.

Let $\HH$ be the periodic Hilbert transform, and set
\[    g=f-\avg f,
    \qquad
    F=-\HH g.
\]
Since $\HH^2g=-g$ for $\avg g=0$,
\begin{equation}
    g=\HH F.
    \label{eq:g-HF}
\end{equation}

For the logarithmic kernel define
$
    \psi_{a,b}:=-\HH\ell_{a,b}.
$
For $0\le a<1$, \eqref{eq:ellab-fourier} gives
\[
    \psi_{a,b}(x)
    =\sum_{m=1}^{\infty}
      \frac{a^m}{m}\sin 2\pi m(x-b),
\]
and hence
\[
    \psi_{a,b}'(x)
    =2\pi\sum_{m=1}^{\infty}a^m\cos 2\pi m(x-b)
    =\pi\bigl(P_a(x-b)-1\bigr).
\]
Since $P_a\ge0$,
\begin{equation}
    \psi_{a,b}'\ge-\pi
    \qquad\text{in distributions.}
    \label{eq:psi-lower-derivative}
\end{equation}
Moreover,
\begin{equation}
    \psi_{a,b}(x)
    =\operatorname{Im}\bigl[-\log(1-ae(x-b))\bigr],
    \qquad
    \norm{\psi_{a,b}}_{L^\infty(\T)}\le\pi.
    \label{eq:psi-bound}
\end{equation}
For $a=1$, the $L^p$ boundary value, $1\le p<\infty$, is the shifted sawtooth
\[
    \psi_{1,b}(x)
    =\pi\Bigl(\frac12-\{x-b\}\Bigr),
    \quad x\ne b,
\]
with
\begin{equation}
    \psi_{1,b}'
    =\pi(\delta_b-1)
    =\pi(P_1(\cdot-b)-1)
    \ge-\pi
    \label{eq:sawtooth-derivative}
\end{equation}
in distributions, where $P_1=\delta_0$.

Let
$
    \phi:=-\HH(h-\avg h).
$
Then
\begin{align}
    \norm{\phi}_{L^\infty(\T)}
    &\le \sum_{k\ne0}|\widehat h(k)|
    \le \sum_{k\ne0}|k|\,|\widehat h(k)|
    \le CK_r,
    \label{eq:phi-bound}\\
    \norm{\phi'}_{L^\infty(\T)}
    &\le 2\pi\sum_{k\ne0}|k|\,|\widehat h(k)|
    \le CK_r.
    \label{eq:phi-prime-bound}
\end{align}
The constants in \eqref{eq:Riesz-used} disappear under $\HH$, and therefore
\[
    F
    =\phi+
      \int\psi_{a_\zeta,b_\zeta}\,d\nu(\zeta).
\]
Using \eqref{eq:psi-bound}, \eqref{eq:phi-bound}, and
$\nu(\T_{r/2})\le Cr^{-1}$,
\[
    \norm{F}_{L^\infty(\T)}\le CK_r.
\]
From \eqref{eq:psi-lower-derivative}, \eqref{eq:sawtooth-derivative}, and
\eqref{eq:phi-prime-bound},
\[
    F'\ge
    -CK_r-\pi\nu(\T_{r/2})
    \ge-CK_r
    \qquad\text{in distributions.}
\]
Choose an absolute constant $A_0>0$ such that
\begin{equation}
    F'\ge-A_0K_r
    \qquad\text{in distributions},
    \label{eq:F-derivative-fixed}
\end{equation}
and set $\Lambda_r:=A_0K_r.$ Then
\[
    (F(x)+\Lambda_rx)'\ge0
\]
on every real lift.  Since a positive distribution is a positive Radon measure,
$F(x)+CK_rx$ admits a nondecreasing representative. Hence, for almost every
$x<y$,
\begin{equation}
    F(y)\ge F(x)-\Lambda_r(y-x).
    \label{eq:one-sided-Lipschitz}
\end{equation}

Choose $\delta\in C_c^\infty((-1,1))$ with
\begin{equation}
    \delta\ge0,
    \qquad
    \int_{\R}\delta=1,
\end{equation}
and for $L\ge2$ let $\delta_L$ be the periodicization of $L\delta(Lx).$ Thus $\delta_L\ge0$, $\int_\T\delta_L=1$, and
\begin{equation}
    \supp\delta_L
    \subset\{x\in\T:\dist(x,\Z)\le L^{-1}\}.
    \label{eq:mollifier-support}
\end{equation}
At every Lebesgue point of $F$, \eqref{eq:one-sided-Lipschitz} gives
\begin{align}
    F(x)
    &\le
    \int_\T
    F(x+L^{-1}-u)\delta_L(u)\,du
    +\frac{2\Lambda_r}{L} \notag\\
    &=(F*\delta_L)(x+L^{-1})+\frac{2\Lambda_r}{L},
    \label{eq:F-upper-mollified}\\
    F(x)
    &\ge
    \int_\T
    F(x-L^{-1}-u)\delta_L(u)\,du
    -\frac{2\Lambda_r}{L} \notag\\
    &=(F*\delta_L)(x-L^{-1})-\frac{2\Lambda_r}{L}.
    \label{eq:F-lower-mollified}
\end{align}
Consequently,
\begin{equation}
    \norm{F}_{L^\infty(\T)}
    \le
    \norm{F*\delta_L}_{L^\infty(\T)}
    +\frac{C K_r}{L}.
    \label{eq:F-from-mollified}
\end{equation}

The periodic Hilbert transform has principal-value kernel $\cot(\pi x)$.
Let $d(x):=\dist(x,\Z).$ If $d(x)\ge2/L$, then \eqref{eq:mollifier-support} gives
\[
    |\HH\delta_L(x)|
    \le
    \int_\T\delta_L(u)|\cot\pi(x-u)|\,du
    \le \frac{C}{d(x)}.
\]
If $d(x)<2/L$ and $L\ge8$, write $x=X/L$, $u=Y/L$ on a suitable lift, with $|X|<2$ and $|Y|<1$.  Since
\[
    \cot\frac{\pi(X-Y)}{L}
    =\frac{L}{\pi(X-Y)}+O\left(\frac{|X-Y|}{L}\right),
\]
and the principal-value Hilbert transform of the fixed smooth function
$\delta$ is bounded on $|X|\le2$, one obtains
\[
    |\HH\delta_L(x)|\le CL.
\]
For $2\le L<8$ the same estimate follows after enlarging the absolute constant. Thus
\[
    |\HH\delta_L(x)|
    \le C\min\{L,d(x)^{-1}\},
\]
and therefore
\begin{align}
    \norm{\HH\delta_L}_{L^\infty(\T)}
    &\le CL,
    \label{eq:Hdelta-Linf}\\
    \norm{\HH\delta_L}_{L^1(\T)}
    &\le
    C\int_0^{2/L}L\,dx
    +C\int_{2/L}^{1/2}x^{-1}\,dx
    \le C\log L.
    \label{eq:Hdelta-L1}
\end{align}

Since $\HH$ commutes with periodic convolution and annihilates constants,
\[
    F*\delta_L
    =-f*(\HH\delta_L).
\]
Using \eqref{eq:bmo-decomp-assumption}, \eqref{eq:Hdelta-Linf}, and
\eqref{eq:Hdelta-L1},
\begin{align}
    \norm{F*\delta_L}_{L^\infty(\T)}
    &\le
    \norm{f_0}_{\infty}\norm{\HH\delta_L}_{1}
    +\norm{f_1}_{1}\norm{\HH\delta_L}_{\infty} \notag\\
    &\le
    C\eps_0\log L+C\eps_1L.
    \label{eq:F-mollified-bound}
\end{align}
Combining \eqref{eq:F-from-mollified} and \eqref{eq:F-mollified-bound},
\[
    \norm{F}_{L^\infty(\T)}
    \le
    C\left(
        \eps_0\log L
        +\eps_1L
        +\frac{K_r}{L}
    \right).
\]
Finally, by \eqref{eq:g-HF} and the periodic bound
$\HH:L^\infty(\T)\to\BMO(\T)$,
\[
    \norm f_{\BMO(\T)}
    =\norm{f-\avg f}_{\BMO(\T)}
    =\norm{\HH F}_{\BMO(\T)}
    \le C\norm{F}_{L^\infty(\T)},
\]
which proves \eqref{eq:bmo-decomp}.

Assume \eqref{eq:l1-smallness}.  If $\eps=0$, then $f=\avg f$ almost
everywhere on $\T$, and \eqref{eq:bmo-l1} is immediate.  Let $\eps>0$.
Subtracting the constant $\avg f$ changes neither the Riesz measure nor the nonzero Fourier coefficients of the harmonic remainder. Hence the above steps can be directly applied to
\[
    g=f-\avg f
\]
with $g_0=0$, $g_1=g$, and give
\begin{equation}
    \norm f_{\BMO(\T)}
    =\norm g_{\BMO(\T)}
    \le C\left(\eps L+\frac{K_r}{L}\right).
    \label{eq:optimization-start}
\end{equation}
Choose
\[
    L=\max\left\{2,\sqrt{\frac{K_r}{\eps}}\right\}.
\]
If $K_r/\eps\ge4$, then
\[
    \eps L+\frac{K_r}{L}
    =2\sqrt{K_r\eps}.
\]
If $K_r/\eps<4$, then $K_r\le4\eps\le4$ and
$\sqrt{K_r\eps}\ge K_r/2$ up to an absolute factor, taking $L=2$ in
\eqref{eq:optimization-start} again yields
\[
    \eps L+\frac{K_r}{L}
    \le C\sqrt{K_r\eps}.
\]
Hence \eqref{eq:bmo-l1} follows.

For a fixed bound $|f|\le B$, apply the argument to $B^{-1}f$ and absorb
$B$ into the constants.
\end{proof}

\begin{remark}
\label{rem:strip-loss}
The estimate \eqref{eq:bmo-decomp} has the form
\begin{equation}
    \norm f_{\BMO}
    \lesssim
    \eps_0\log L+\eps_1L+K_rL^{-1}.
\end{equation}
In particular, the strip cost $K_r$ is attached only to the last term; it
does not multiply $\eps_0\log L$.  This separation is the input needed in
the shrinking-strip boosting argument of Section~3.
\end{remark}

% Section 3
\section{Several-variable boosting on shrinking strips}
\label{sec:boosting}

For separately subharmonic functions, positive averaging in real variables preserves the
complex strips of all retained coordinates.  We record the precise form used below.
Let
\[
 \Omega:=\prod_{j=1}^{m}\Tstrip{r_j},
 \qquad 0<r_j<\frac1{10}.
\]
Fix $J\subset\{1,\dots,m\}$ and a probability measure $\pi$ on $\T^J$.
For $t=(t_j)_{j\in J}\in\T^J$, let $t_J\in\T^m$ denote the vector whose
$j$th coordinate is $t_j$ for $j\in J$ and $0$ otherwise. And let
\[
 \Omega_{J^c}:=\prod_{j\notin J}\Tstrip{r_j}.
\]
We prove the following lemma.
\begin{lemma}\label{lem:positive-average}
Let $B>0$ and $f:\Omega\to\R$ be continuous and separately subharmonic, with
\[
 |f(z)|\le B,\qquad z\in\Omega.
\]
Define
\begin{equation}\label{eq:CJpi}
 (\CJpi f)(z)
 :=\int_{\T^J} f(z-t_J)\,d\pi(t),
 \qquad z\in\Omega.
\end{equation}
If $\iota_J(t,z_{J^c})\in\Omega$ denotes the point whose coordinates in $J$
are $t$ and whose remaining coordinates are $z_{J^c}$, define
\begin{equation}\label{eq:AJpi}
 (\AJpi f)(z_{J^c})
 :=\int_{\T^J} f\bigl(\iota_J(t,z_{J^c})\bigr)\,d\pi(t),
 \qquad z_{J^c}\in\Omega_{J^c}.
\end{equation}
Then $\CJpi f$ is continuous and separately subharmonic on $\Omega$, while
$\AJpi f$ is continuous and separately subharmonic on $\Omega_{J^c}$.  Moreover,
\begin{equation}\label{eq:average-bound}
 \|\CJpi f\|_{L^\infty(\Omega)}\le B,
 \qquad
 \|\AJpi f\|_{L^\infty(\Omega_{J^c})}\le B.
\end{equation}
For every retained coordinate $j$, the one-dimensional sections have the same
strip width $r_j$ and satisfy the estimates of Lemmas \ref{lem:riesz}-\ref{lem:oneD} with strip cost
\begin{equation}\label{eq:coordinate-cost}
 K_j\le C B r_j^{-4}.
\end{equation}
The same conclusion holds after any finite succession of the operations
\eqref{eq:CJpi}--\eqref{eq:AJpi}; no additional loss depending on the number of
operations is introduced.
\end{lemma}

\begin{proof}
Write $g$ for either $\CJpi f$ or $\AJpi f$, and call a coordinate $j$
\emph{retained} if it belongs to the domain of $g$.  For fixed $t\in\T^J$,
let $\Phi_t$ denote respectively
\[
 \Phi_t z=z-t_J
 \qquad\text{or}\qquad
 \Phi_t z_{J^c}=\iota_J(t,z_{J^c}).
\]
For every retained $j$ and every admissible $\xi\in\C$,
\begin{equation}\label{eq:Phi-coordinate}
 \Phi_t(z+\xi e_j)=\Phi_t z+\xi e_j,
 \qquad
 \Im(\Phi_t z)_j=\Im z_j.
\end{equation}
Hence the $j$th strip $\Tstrip{r_j}$ is unchanged.  Since $\pi$ is a
probability measure,
\begin{equation}\label{eq:sup-preserved}
 |g(z)|
 \le \int_{\T^J}|f(\Phi_t z)|\,d\pi(t)
 \le B.
\end{equation}
This proves \eqref{eq:average-bound}.

Fix a retained coordinate $j$, all other coordinates, and a disk $ D(z_j,a)\Subset\Tstrip{r_j}.$ By \eqref{eq:Phi-coordinate},
\[
 D\bigl((\Phi_t z)_j,a\bigr)\Subset\Tstrip{r_j},
 \qquad t\in\T^J.
\]
Separate subharmonicity of $f$ gives
\begin{align}
 f(\Phi_t z)
 &\le
 \frac1{2\pi}\int_0^{2\pi}
 f\bigl(\Phi_t z+a e^{i\theta}e_j\bigr)\,d\theta .
 \label{eq:submean-before-average}
\end{align}
Integrating \eqref{eq:submean-before-average} against $d\pi(t)$ and using
Fubini's theorem,
\begin{align}
 g(z)
 &\le
 \frac1{2\pi}\int_0^{2\pi}
 \left[
 \int_{\T^J}
 f\bigl(\Phi_t(z+a e^{i\theta}e_j)\bigr)\,d\pi(t)
 \right]d\theta \notag\\
 &=
 \frac1{2\pi}\int_0^{2\pi}
 g(z+a e^{i\theta}e_j)\,d\theta .
 \label{eq:submean-after-average}
\end{align}
Thus $g$ is subharmonic in every retained complex coordinate.

We next verify continuity.  Let $K$ be a compact subset of the domain of $g$ and set
\[
 \widetilde K
 :=\{\Phi_t z:\ z\in K,\ t\in\T^J\}\subset\Omega.
\]
The map $(z,t)\mapsto\Phi_t z$ is continuous, $K\times\T^J$ is compact, and
therefore $\widetilde K\Subset\Omega$.  Hence $f$ is uniformly continuous on
$\widetilde K$.  For $z,w\in K$, when $w\to z$,
\begin{align*}
 |g(z)-g(w)|
 &\le
 \int_{\T^J}|f(\Phi_t z)-f(\Phi_t w)|\,d\pi(t) \notag\\
 &\le
 \sup_{t\in\T^J}|f(\Phi_t z)-f(\Phi_t w)|
 \longrightarrow0
\end{align*}
Thus $g$ is continuous.

Fix a retained $j$ and freeze all other variables.  The resulting section
\[
 u_j(\zeta)
 :=g(z_1,\dots,z_{j-1},\zeta,z_{j+1},\dots)
\]
is continuous and subharmonic on $\Tstrip{r_j}$, with $|u_j|\le B$. Applying Lemmas \ref{lem:riesz}-\ref{lem:oneD} to $B^{-1}u_j$ and rescaling gives
\begin{equation}\label{eq:section-cost}
 K_j=C_*B r_j^{-4},
\end{equation}
where $C_*$ is absolute.  In particular, the one-sided Poisson estimate reads
\begin{equation}\label{eq:section-poisson}
 u_j(x)\le P_tu_j(x)+C_*B r_j^{-4}t,
 \qquad 0<t<1,
\end{equation}
and, whenever $u_j=u_{j,0}+u_{j,1}$ on $\T$ with
$
 \|u_{j,0}\|_\infty\le\varepsilon_0,
 \quad
 \|u_{j,1}\|_1\le\varepsilon_1,
$
the corresponding BMO estimate is
\[
 \|u_j\|_{\BMO(\T)}
 \le C\left(
 \varepsilon_0\log L+\varepsilon_1L+
 \frac{C_*B r_j^{-4}}{L}
 \right),
 \qquad L\ge2.
\]

After either operation, the retained widths $r_j$ and the bound $B$ remain
unchanged by \eqref{eq:Phi-coordinate} and \eqref{eq:sup-preserved}.  Repeating
\eqref{eq:section-cost} after each subsequent operation therefore gives the same
constant $C_*B r_j^{-4}$, independently of the number of averaging steps.
\end{proof}

\subsection{Coordinatewise smoothing}
\begin{lemma}\label{lem:fejer-remainder}
Let $m\ge 1$, $\rho>0$, and let $f$ be continuous and separately subharmonic on
\[
\T_{\rho/N}^{m}=\prod_{j=1}^{m}\T_{\rho/N},
\qquad |f|\le 1.
\]
Write $x=(x',t)\in\T^{m-1}\times\T$. Let $F_K$ be the Fej\'er kernel (see \cite{SteinShakarchi,Zygmund}),
\[
\widehat F_K(k)=\left(1-\frac{|k|}{K}\right)_{+},
\qquad \int_{\T}F_K=1,
\qquad K\in\Z,\ K\ge2,
\]
and, for $1\le j\le m-1$, let
$
\mathcal F_j g:=F_K*_{x_j}g,
\quad
v_K:=\mathcal F_1\cdots\mathcal F_{m-1}f.
$
There exist $C=C(m,\rho)$ and $N_0=N_0(m,\rho)$ such that, for every $N\ge N_0$,
\begin{equation}\label{eq:fejer-L1-sharp}
\sup_{t\in\T}
\norm{f(\cdot,t)-v_K(\cdot,t)}_{L^1(\T^{m-1})}
\le C N^{4}K^{-1/2}.
\end{equation}
Consequently, after increasing $N_0$, one may use the fixed exponent $p=6$ and write
\begin{equation}\label{eq:fejer-L1-p6}
\sup_{t\in\T}
\norm{f(\cdot,t)-v_K(\cdot,t)}_{L^1(\T^{m-1})}
\le N^{p}K^{-1/2}.
\end{equation}
Moreover,
\begin{equation}\label{eq:fejer-gradient}
\norm{\nabla_{x'}v_K}_{L^\infty(\T^m)}
\le C_m K^2.
\end{equation}
For $m=1$, the assertions are understood with $x'$ empty and $v_K=f$.
\end{lemma}

\begin{proof}
Assume $m\ge2$ and fix $1\le j\le m-1$. For
\[
x_{\widehat j}
=(x_1,\ldots,x_{j-1},x_{j+1},\ldots,x_{m-1})\in\T^{m-2},
\qquad t\in\T,
\]
define the partial Fourier coefficient
\[
\widehat f_j(k;x_{\widehat j},t)
:=
\int_{\T}
 f(x',t)e(-kx_j)\,dx_j,
\qquad k\in\Z.
\]
For fixed $(x_{\widehat j},t)$, the map
\[
z_j\longmapsto
f(x_1,\ldots,x_{j-1},z_j,x_{j+1},\ldots,x_{m-1},t)
\]
is subharmonic on $\T_{\rho/N}$ and bounded by one. Lemma~\ref{lem:riesz} therefore gives, for $k\ne0$,
\[
\bigl|
\widehat f_j(k;x_{\widehat j},t)
\bigr|
\le
C\left(\frac{\rho}{N}\right)^{-4}|k|^{-1}
\le
C_{\rho}N^4|k|^{-1},
\]
uniformly in $(x_{\widehat j},t)$.

Set $f_t(x'):=f(x',t)$. Since
\[
1-\widehat F_K(k)
=
\begin{cases}
|k|/K,&1\le |k|<K,\\
1,&|k|\ge K,\\
0,&k=0,
\end{cases}
\]
Parseval in the $x_j$-variable gives, for fixed $(x_{\widehat j},t)$,
\begin{align}
\norm{(I-\mathcal F_j)f_t
}_{L^2_{x_j}}^2
&=\sum_{k\in\Z}
|1-\widehat F_K(k)|^2
\bigl|
\widehat f_j(k;x_{\widehat j},t)
\bigr|^2
\notag\\
&\le
C_{\rho}N^8
\left[
\sum_{1\le |k|<K}
\frac{|k|^2}{K^2}\frac1{k^2}
+
\sum_{|k|\ge K}
\frac1{k^2}
\right]
\notag\\
&\le
C_{\rho}N^8
\left(
\frac{K}{K^2}
+
\frac1K
\right)
\le
C_{\rho}\frac{N^8}{K}.
\label{eq:one-coordinate-L2}
\end{align}
Integrating in $x_{\widehat j}$ yields
\begin{equation}\label{eq:one-coordinate-L2-global}
\norm{(I-\mathcal F_j)f_t}_{L^2(\T^{m-1})}
\le
C_{\rho}N^4K^{-1/2},
\end{equation}
uniformly in $t$ and $j$.

Since each $F_K$ is nonnegative and has integral one,
\begin{equation}\label{eq:L2-contraction}
\norm{\mathcal F_j g}_{L^2(\T^{m-1})}
\le
\norm{g}_{L^2(\T^{m-1})}.
\end{equation}
Moreover,
\[
I-\mathcal F_1\cdots\mathcal F_{m-1}
=
\sum_{j=1}^{m-1}
\mathcal F_1\cdots\mathcal F_{j-1}(I-\mathcal F_j).
\]
Hence, by \eqref{eq:L2-contraction} and \eqref{eq:one-coordinate-L2-global},
\[
\norm{f_t-v_K(\cdot,t)}_{L^2(\T^{m-1})}
\le
\sum_{j=1}^{m-1}
\norm{(I-\mathcal F_j)f_t}_{L^2(\T^{m-1})}\le
C_{m,\rho}N^4K^{-1/2}.
\]
Since $\T^{m-1}$ has normalized measure one,
\[
\norm{g}_{L^1(\T^{m-1})}
\le
\norm{g}_{L^2(\T^{m-1})},
\]
and therefore
\[
\sup_{t\in\T}
\norm{f(\cdot,t)-v_K(\cdot,t)}_{L^1(\T^{m-1})}
\le
C_{m,\rho}N^4K^{-1/2}.
\]
This proves \eqref{eq:fejer-L1-sharp}. Choose $N_0=N_0(m,\rho)$ so that
\[
C_{m,\rho}N^4\le N^6
\qquad (N\ge N_0).
\]
Then \eqref{eq:fejer-L1-p6} follows with $p=6$.

For $1\le j\le m-1$,
\begin{equation}\label{eq:derivative-convolution}
\partial_{x_j}v_K
=
(F_K')*_{x_j}
\prod_{\substack{1\le i\le m-1\\ i\ne j}}
(F_K*_{x_i})f.
\end{equation}
By the definition of $F_K$,
\[
F_K(x)
=
\sum_{|k|<K}
\left(1-\frac{|k|}{K}\right)e(kx),
\]
we have
\[
\norm{F_K'}_{L^1(\T)}
\le
\norm{F_K'}_{L^\infty(\T)}
\le
2\pi
\sum_{|k|<K}
|k|\left(1-\frac{|k|}{K}\right)
\le
CK^2.
\]
Also
\[
\norm{F_K}_{L^1(\T)}=1,
\qquad
\norm{f}_{L^\infty(\T^m)}\le1.
\]
Young's inequality and \eqref{eq:derivative-convolution} give
\[
\norm{\partial_{x_j}v_K}_{L^\infty(\T^m)}
\le
\norm{F_K'}_1
\prod_{i\ne j}\norm{F_K}_1
\norm{f}_\infty
\le
CK^2.
\]
Summing over $1\le j\le m-1$ yields \eqref{eq:fejer-gradient}.
\end{proof}

\begin{proposition}\label{prop:boosting}
Fix $m\ge1$, $s>0$, and $\rho>0$. There exist a finite number $P_m(s)$,
constants $C,c>0$, and $N_*=N_*(m,s,\rho)$ such that the following holds
for every integer $N\ge N_*$: Let $f$ be continuous and separately
subharmonic on $\T_{\rho/N}^{m}$, and $|f|\le1$,
and suppose that, on the real torus,
\begin{equation}\label{eq:boost-decomp}
f=f_0+f_1,
\qquad
\norm{f_0}_{L^\infty(\T^m)}\le N^{-s},
\qquad
\norm{f_1}_{L^1(\T^m)}\le N^{-P_m(s)}.
\end{equation}
Then
\begin{equation}\label{eq:boost-conclusion}
\mes\left\{x\in\T^m:
|f(x)|>C N^{-\beta_m(s)}\right\}
\le
C\exp\bigl(-cN^{\beta_m(s)}\bigr),
\end{equation}
where one may take
\begin{equation}\label{eq:beta-def}
\beta_m(s)=\frac{s}{4\cdot 8^{m-1}}.
\end{equation}
Neither $f_0$ nor $f_1$ is required to be subharmonic.
\end{proposition}

\begin{proof}
Fix the exponent $p=6$ from the uniform Fej\'er remainder. After increasing
$N_*$, the one-dimensional strip cost at width $\rho/N$ satisfies
\begin{equation}\label{eq:strip-cost-p}
K_{\rho/N}\le N^p.
\end{equation}
Set
\begin{equation}\label{eq:P1-def}
P_1(s):=p+2s+6.
\end{equation}
We proceed by induction on $m$. Let $m=1$. Apply the BMO estimate with
\[
\varepsilon_0=N^{-s},
\qquad
\varepsilon_1=N^{-P_1(s)},
\qquad
L=N^{p+s+2}.
\]
Using \eqref{eq:strip-cost-p} and \eqref{eq:P1-def},
\begin{align}
\norm{f}_{\mathrm{BMO}(\T)}
&\le
C\left(
N^{-s}\log L
+N^{-P_1(s)}L
+\frac{N^p}{L}
\right)
\notag\\
&\le
C\left(
N^{-s}\log N
+N^{-s-4}
+N^{-s-2}
\right)
\notag\\
&\le C N^{-s/2}.
\label{eq:base-bmo}
\end{align}
Moreover,
\[
|\langle f\rangle|
\le
\norm{f_0}_\infty+\norm{f_1}_1
\le
N^{-s}+N^{-P_1(s)}
\le 2N^{-s}.
\]
For sufficiently large $N$,
\[
2N^{-s}\le \frac12 N^{-s/4}.
\]
Hence
\begin{align}
\mes\{|f|>2N^{-s/4}\}
&\le
\mes\left\{|f-\langle f\rangle|>N^{-s/4}\right\}
\notag\\
&\le
C\exp\left(
-c\frac{N^{-s/4}}{N^{-s/2}}
\right)
\le
C e^{-cN^{s/4}}.
\label{eq:boost-base-JN}
\end{align}
Thus \eqref{eq:boost-conclusion} holds with $\beta_1(s)=\frac{s}{4}.$

Assume the proposition in dimension $m-1$. Put $s':=s/8,$ and choose, in this order,
\[
q:=\left\lceil
2\bigl(p+P_{m-1}(s')+3\bigr)
\right\rceil,
\]
\[
P_m(s):=q(m-1)+P_1(s)+2.
\]
Let
\[
K:=\lceil N^q\rceil,
\qquad
\F_jg:=F_K*_{x_j}g,
\qquad
v_K:=\F_1\cdots\F_{m-1}f.
\]
On the real torus define
\begin{equation}\label{eq:vk-components}
v_{K,\ell}:=\F_1\cdots\F_{m-1}f_\ell,
\qquad \ell=0,1.
\end{equation}
Then
$
v_K=v_{K,0}+v_{K,1}.
$
Positivity and unit mass of $F_K$ give
\[
\norm{v_{K,0}}_\infty
\le
\norm{f_0}_\infty
\le N^{-s}.
\]
Since
\[
0\le F_K\le K,
\qquad
\norm{F_K^{\otimes(m-1)}}_\infty=K^{m-1},
\]
Fubini's theorem gives, for every $x'\in\T^{m-1}$,
\begin{align}
\norm{v_{K,1}(x',\cdot)}_{L^1_t}
&\le
K^{m-1}\norm{f_1}_{L^1(\T^m)}
\notag\\
&\le
K^{m-1}N^{-P_m(s)}
\notag\\
&\le
2^{m-1}N^{q(m-1)-P_m(s)}
\notag\\
&=
2^{m-1}N^{-P_1(s)-2}
\le
N^{-P_1(s)}
\label{eq:vk1-fiber}
\end{align}
for all sufficiently large $N$.

By positive real convolution, $v_K$ is separately subharmonic on
$\T_{\rho/N}^m$ and $|v_K|\le1$. Applying the case \(m=1\) to the $t$-variable,
uniformly in $x'$, yields
\begin{equation}\label{eq:fiber-LDT}
\sup_{x'\in\T^{m-1}}
\mes_t\left\{
|v_K(x',t)|>C_0N^{-s/4}
\right\}
\le
C_0e^{-c_0N^{s/4}}.
\end{equation}

The Fej\'er derivative bound gives
\begin{equation}\label{eq:vk-gradient}
\norm{\nabla_{x'}v_K}_\infty\le C_mK^2.
\end{equation}
Choose
$
\eta:=c_mK^{-2}N^{-s/4},
$
with $c_m>0$ so that
\begin{equation}\label{eq:lip-error}
C_mK^2\sqrt{m-1}\,\eta
\le N^{-s/4}.
\end{equation}
Let $\mathcal N\subset\T^{m-1}$ be a rectangular $\eta$-net. Then
\begin{equation}\label{eq:net-card}
|\mathcal N|
\le C_m\eta^{-(m-1)}
\le
C_m\bigl(K^2N^{s/4}\bigr)^{m-1}
\le
C_mN^{D_q},
\end{equation}
where $D_q:=\left(2q+\frac{s}{4}\right)(m-1)$ and $K\le2N^q$ has been used.
Define
\[
E_N
:=
\bigcup_{\xi\in\mathcal N}
\left\{
 t\in\T:
 |v_K(\xi,t)|>C_0N^{-s/4}
\right\}.
\]
By \eqref{eq:fiber-LDT} and \eqref{eq:net-card},
\begin{align}
\mes(E_N)
&\le
C_mN^{D_q}e^{-c_0N^{s/4}}
\notag\\
&\le
C_m\exp\left(D_q\log N-c_0N^{s/4}\right)
\notag\\
&\le
C_me^{-c_1N^{s/4}}
\label{eq:EN-measure}
\end{align}
for sufficiently large $N$.
For $t\notin E_N$ and $x'\in\T^{m-1}$ choose
$\xi\in\mathcal N$ with $|x'-\xi|\le\sqrt{m-1}\,\eta$.
Then \eqref{eq:vk-gradient}--\eqref{eq:lip-error} give
\[
|v_K(x',t)|
\le
|v_K(\xi,t)|
+C_mK^2|x'-\xi|
\le
(C_0+1)N^{-s/4}.
\]
After increasing $N_*$,
\[
(C_0+1)N^{-s/4}
\le
N^{-s/8}
=N^{-s'}.
\]

For every $t\in\T$, the uniform Fej\'er remainder gives
\begin{equation}
\norm{f(\cdot,t)-v_K(\cdot,t)}_{L^1(\T^{m-1})}\le
N^pK^{-1/2}\le
N^{-P_{m-1}(s')-3}
\le
N^{-P_{m-1}(s')}.
\label{eq:remainder-induction}
\end{equation}
For $t\notin E_N$, the section $x'\longmapsto f(x',t)$ is continuous and separately subharmonic on
$\T_{\rho/N}^{m-1}$, with $|f|\le1$, and therefore the induction hypothesis gives
\begin{equation}\label{eq:section-boost}
\mes_{x'}\left\{
|f(x',t)|>C N^{-\beta_{m-1}(s')}
\right\}
\le
C\exp\left(-cN^{\beta_{m-1}(s')}\right)
\end{equation}
uniformly for $t\notin E_N$.
Since
\[
\beta_{m-1}(s')
=
\frac{s/8}{4\cdot 8^{m-2}}
=
\frac{s}{4\cdot 8^{m-1}}
=
\beta_m(s)
<\frac{s}{4},
\]
Fubini's theorem, \eqref{eq:EN-measure}, and \eqref{eq:section-boost} yield
\begin{align}
&\mes\left\{
(x',t)\in\T^m:
|f(x',t)|>C N^{-\beta_m(s)}
\right\}
\notag\\
&\qquad\le
\mes(E_N)
+
C\exp\left(-cN^{\beta_m(s)}\right)
\notag\\
&\qquad\le
C\exp\left(-cN^{\beta_m(s)}\right).
\end{align}
This proves \eqref{eq:boost-conclusion}.

All parameter choices are made in the order $P_{m-1}(s/8)$, then $q$, then $P_m(s)$, so no parameter is defined circularly.
\end{proof}

\begin{lemma}\label{lem:one-sided-concentration}
Fix $m\ge1$ and $s>0$. Let $f$ be a continuous real-valued separately subharmonic
function on a product complex neighborhood of $\T^m$, with a fixed absolute bound there.
For $J\subset\{1,\dots,m\}$, write
\[
(\mathcal A_Jf)(x_{J^c})
:=\int_{\T^J}f(t_J,x_{J^c})\,dt_J .
\]
Assume that there is a quantity $\mathcal K_N$ such that
\begin{equation}\label{eq:K-growth}
1\le \mathcal K_N\le C_1N^\kappa,
\end{equation}
and that every one-dimensional real section $u$ of every $\mathcal A_Jf$
(including $J=\varnothing$) satisfies the $L^1$-to-BMO estimate
\begin{equation}\label{eq:section-BMO}
\|u\|_{\BMO(\T)}
\le C_1\bigl(\mathcal K_N
\|u-\avg{u}\|_{L^1(\T)}\bigr)^{1/2}
\end{equation}
whenever $\|u-\avg{u}\|_{L^1(\T)}\le1$.
Suppose, on the real torus,
\begin{equation}\label{eq:one-sided-hyp}
f\le C_0N^{-s},
\qquad
|\avg f|\le C_0N^{-s},
\qquad
\kappa\le \frac{s}{8^{m+1}}.
\end{equation}
Then there are $C,c>0$ and $N_*=N_*(m,s,C_0,C_1)>0$ such that,
for $N\ge N_*$,
\begin{equation}\label{eq:one-sided-tail}
\mes\bigl\{x\in\T^m:\ |f(x)|>CN^{-s/8^m}\bigr\}
\le C\exp\bigl(-cN^{s/8^m}\bigr).
\end{equation}

If $f$ is bounded on a product strip and every retained coordinate has
strip cost at most $C_1N^\kappa$, then \eqref{eq:section-BMO} for all
intermediate averages follows from Lemmas \ref{lem:oneD} and \ref{lem:positive-average}.
\end{lemma}

\begin{proof}
We argue by induction on $m$. Let $m=1$ and set
$
g:=f-\avg f,
\quad
\avg g=0.
$
Since $\int_\T g=0$,
\[
\|g\|_{L^1(\T)}
=2\int_\T g_+\,dx.
\]
By \eqref{eq:one-sided-hyp},
\[
g(x)
=f(x)-\avg f
\le C_0N^{-s}+|\avg f|
\le 2C_0N^{-s},
\]
hence
\begin{equation}\label{eq:base-L1}
\|f-\avg f\|_{L^1(\T)}
\le 4C_0N^{-s}.
\end{equation}
For $N$ large, \eqref{eq:section-BMO}, \eqref{eq:K-growth} and
\eqref{eq:base-L1} give
\begin{equation}
\|f\|_{\BMO(\T)}
\le C\bigl(N^\kappa N^{-s}\bigr)^{1/2}
 = C N^{-(s-\kappa)/2}\le C N^{-s/4},
\label{eq:base-BMO}
\end{equation}
where $\kappa\le s/8^2<s/2$.
John--Nirenberg inequality, with $u=N^{-s/8}$, yields
\[
\mes\bigl\{|f-\avg f|>N^{-s/8}\bigr\}
\le C\exp\left(
-c\frac{N^{-s/8}}{N^{-s/4}}
\right)\le C\exp\bigl(-cN^{s/8}\bigr).
\]
Since
$
|\avg f|\le C_0N^{-s}=o(N^{-s/8}),
$
\eqref{eq:one-sided-tail} follows for $m=1$ after enlarging the
constant in the threshold.

Assume the result in dimension $m-1$, with $m\ge2$, and define
\begin{equation}\label{eq:g-average}
g(x'):=\int_\T f(x_1,x')\,dx_1,
\qquad
x'=(x_2,\dots,x_m)\in\T^{m-1}.
\end{equation}
Then
$
g\le C_0N^{-s},
\quad
|\avg g|=|\avg f|\le C_0N^{-s}.
$
The section-cost hypothesis is stable under this averaging, and
\[
\kappa\le \frac{s}{8^{m+1}}
<\frac{s}{8^m}.
\]
The induction hypothesis therefore gives, with $t_m:=\frac{s}{8^{m-1}}$, a set $E\subset\T^{m-1}$ such that
\begin{equation}\label{eq:base-bad-set}
\mes(E)\le C\exp(-cN^{t_m}),
\qquad
|g(x')|\le C N^{-t_m},
\quad x'\notin E.
\end{equation}

Fix $x'\notin E$ and set
$
u_{x'}(x_1):=f(x_1,x').
$
Then $\avg{u_{x'}}=g(x').$ Since $\int_\T(u_{x'}-g(x'))=0$,
\[
\|u_{x'}-g(x')\|_{L^1(\T)}
=2\int_\T (u_{x'}-g(x'))_+\,dx_1.
\]
Using \eqref{eq:one-sided-hyp}, \eqref{eq:base-bad-set}, and
$t_m\le s$,
\[
u_{x'}(x_1)-g(x')
\le C_0N^{-s}+|g(x')|
\le C N^{-t_m}.
\]
Hence
\[
\|u_{x'}-g(x')\|_{L^1(\T)}
\le C N^{-t_m}.
\]
For $N$ large, \eqref{eq:section-BMO} and \eqref{eq:K-growth} give
\[
\|u_{x'}\|_{\BMO(\T)}
\le C\bigl(N^\kappa N^{-t_m}\bigr)^{1/2}
= C N^{-(t_m-\kappa)/2}.
\]
Since
\[
\kappa\le \frac{s}{8^{m+1}}
=\frac{t_m}{64}
<\frac{t_m}{2},
\]
we obtain
\[
\|u_{x'}\|_{\BMO(\T)}
\le C N^{-t_m/4}.
\]
John--Nirenberg inequality at the threshold $N^{-t_m/8}$ yields, uniformly for $x'\notin E$,
\[
\mes_{x_1}\Bigl\{
|u_{x'}-g(x')|>N^{-t_m/8}
\Bigr\}\le C\exp\left(
-c\frac{N^{-t_m/8}}{N^{-t_m/4}}
\right)\le C\exp\bigl(-cN^{t_m/8}\bigr).
\]
Moreover,
\[
|g(x')|\le C N^{-t_m}=o(N^{-t_m/8}).
\]
Thus, after enlarging the threshold constant,
\begin{equation}\label{eq:fiber-final}
\mes_{x_1}\Bigl\{
|f(x_1,x')|>C N^{-t_m/8}
\Bigr\}
\le C\exp\bigl(-cN^{t_m/8}\bigr),
\qquad x'\notin E.
\end{equation}

By \eqref{eq:base-bad-set}, \eqref{eq:fiber-final}, and Fubini's theorem,
\begin{align}
&\mes\Bigl\{x\in\T^m:
|f(x)|>C N^{-t_m/8}\Bigr\}
\notag\\
&\qquad\le
\mes(E)
+
\int_{\T^{m-1}\setminus E}
\mes_{x_1}\Bigl\{
|f(x_1,x')|>C N^{-t_m/8}
\Bigr\}\,dx'
\notag\\
&\qquad\le
C e^{-cN^{t_m}}
+C e^{-cN^{t_m/8}}
\le
C e^{-cN^{t_m/8}}.
\label{eq:Fubini-close}
\end{align}
Finally, $\frac{t_m}{8}=\frac{s}{8^m}$, which proves \eqref{eq:one-sided-tail}.
\end{proof}

% Section 4
\section{Quantitative skew-shift averaging}
\label{sec:skew-averaging}

\subsection{Vector Weyl bounds and ergodic preliminaries}
\label{subsec:weyl-ergodic}

We begin with the arithmetic estimate that controls individual Fourier
modes along the skew-shift orbit. Because the orbit phase is quadratic
in the time variable, a single van der Corput differencing step reduces
the relevant exponential sum to a linear phase, where the vector
Diophantine condition can be applied directly. This estimate is the
Fourier input for the short-block large-deviation argument below.

\begin{lemma}\label{lem:vector-weyl}
Assume
\[
\omega\in \DC_d(\gamma,\tau),
\qquad
\distT{k\cdot\omega}\ge \gamma |k|^{-\tau},
\quad 0\ne k\in\Z^d.
\]
Set
\[
\delta:=\frac{1}{2(\tau+1)}.
\]
Then, for every $(k,\ell)\in\Z^d\times\Z^d\setminus\{(0,0)\}$ and every integer $R\ge1$,
\begin{equation}\label{eq:vector-weyl-main}
\sup_{x,y\in\T^d}
\left|
\frac1R\sum_{j=0}^{R-1}
\e^{2\pi i ((k,\ell)\cdot T_\omega^j(x,y))}
\right|
\le
C_{\gamma,\tau}
(1+|k|+|\ell|)^\tau R^{-\delta}.
\end{equation}
Here $|\cdot|$ denotes the $\ell^1$ norm on $\Z^d$.
\end{lemma}

\begin{proof}
For $j\ge0$,
\begin{equation}\label{eq:orbit-iterate}
T_\omega^j(x,y)
=
\left(
 x+j\omega,
 y+jx+\frac{j(j-1)}2\omega
\right).
\end{equation}
Hence
\begin{equation}
(k,\ell)\cdot T_\omega^j(x,y)=k\cdot x+\ell\cdot y+Q(j),\label{eq:phase-split}
\end{equation}
where
\begin{equation}
Q(j)=\frac12(\ell\cdot\omega)j(j-1)
+j(\ell\cdot x+k\cdot\omega).
\label{eq:quadratic-phase}
\end{equation}
The constant phase in \eqref{eq:phase-split} is irrelevant for the modulus.

If $\ell=0$, then $k\ne0$. By the finite geometric-sum identity,
\[
\sum_{j=0}^{R-1}e^{2\pi i j\theta}
=
e^{\pi i(R-1)\theta}
\frac{\sin(\pi R\theta)}{\sin(\pi\theta)},
\qquad
\theta\notin\mathbb Z.
\]
Using
\[
|\sin(\pi\theta)|
\ge
2\|\theta\|_{\mathbb T},
\qquad
\|\theta\|_{\mathbb T}\le \frac12,
\]
we obtain
\begin{equation}
\left|
\frac1R\sum_{j=0}^{R-1}\e^{2\pi i j k\cdot\omega}
\right|\le
\min\left\{1,\frac{C}{R\distT{k\cdot\omega}}\right\}\le
C\gamma^{-1}|k|^\tau R^{-1}\le
C_{\gamma,\tau}(1+|k|)^\tau R^{-\delta}.
\label{eq:linear-case}
\end{equation}
It remains to consider $\ell\ne0$.

Put $a_j:=\e^{2\pi iQ(j)}, 0\le j<R$, and extend $a_j$ by zero to $j\in\Z$. For $1\le H\le R$,
\begin{equation}\label{eq:vdC-identity}
H\sum_{j\in\Z}a_j
=
\sum_{j\in\Z}\sum_{h=0}^{H-1}a_{j+h}.
\end{equation}
The outer sum in \eqref{eq:vdC-identity} has at most $R+H-1\le2R$ nonzero terms. Therefore
\begin{align}
H^2\left|\sum_j a_j\right|^2
&\le
2R\sum_j
\left|\sum_{h=0}^{H-1}a_{j+h}\right|^2\notag\\
&=2R\left[
H\sum_j|a_j|^2
+2\operatorname{Re}\sum_{h=1}^{H-1}(H-h)
\sum_j a_{j+h}\overline{a_j}
\right].
\label{eq:vdC-expand}
\end{align}
Since $\sum_j|a_j|^2=R$,
\begin{equation}\label{eq:vdC-normalized}
\left|
\frac1R\sum_{j=0}^{R-1}a_j
\right|^2
\le
\frac2H
+
\frac4{HR}
\sum_{h=1}^{H-1}
\left|
\sum_{j=0}^{R-h-1}
a_{j+h}\overline{a_j}
\right|.
\end{equation}

By \eqref{eq:quadratic-phase},
\begin{align}
Q(j+h)-Q(j)
&=
 h(\ell\cdot\omega)j
+\frac{h(h-1)}2\,\ell\cdot\omega
+h(\ell\cdot x+k\cdot\omega).
\label{eq:difference-phase}
\end{align}
Thus the $j$-coefficient is $h\ell\cdot\omega$, and geometric summation gives
\begin{align}
\left|
\sum_{j=0}^{R-h-1}a_{j+h}\overline{a_j}
\right|
&\le
\min\left\{R,\frac{C}{\distT{h\ell\cdot\omega}}\right\}.
\label{eq:correlation-geometric}
\end{align}
Since $h\ell\in\Z^d\setminus\{0\}$,
\begin{equation}\label{eq:small-divisor}
\distT{h\ell\cdot\omega}
\ge
\gamma |h\ell|^{-\tau}
=
\gamma h^{-\tau}|\ell|^{-\tau}.
\end{equation}
Substituting \eqref{eq:correlation-geometric}--\eqref{eq:small-divisor} into \eqref{eq:vdC-normalized},
\begin{align}
\left|
\frac1R\sum_{j=0}^{R-1}a_j
\right|^2
&\le
C\left[
H^{-1}
+
\frac1H\sum_{h=1}^{H-1}
\min\left\{1,
\frac{\gamma^{-1}|\ell|^\tau h^\tau}{R}
\right\}
\right]\notag\\
&\le
C\left(
H^{-1}
+
\frac{\gamma^{-1}|\ell|^\tau}{HR}
\sum_{h=1}^{H-1}h^\tau
\right)\notag\\
&\le
C_{\gamma,\tau}
\left(
H^{-1}
+
|\ell|^\tau\frac{H^\tau}{R}
\right).
\label{eq:vdC-diophantine}
\end{align}

For $R\ge2$, take
\[
H:=\left\lfloor R^{1/(\tau+1)}\right\rfloor.
\]
Then
\[
H^{-1}\le C R^{-1/(\tau+1)},
\qquad
\frac{H^\tau}{R}
\le R^{-1/(\tau+1)}.
\]
Hence \eqref{eq:vdC-diophantine} yields
\begin{align}
\left|
\frac1R\sum_{j=0}^{R-1}a_j
\right|^2
&\le
C_{\gamma,\tau}(1+|\ell|^\tau)
R^{-1/(\tau+1)},\notag\\
\left|
\frac1R\sum_{j=0}^{R-1}a_j
\right|
&\le
C_{\gamma,\tau}
(1+|\ell|)^{\tau/2}
R^{-1/(2(\tau+1))}.
\label{eq:quadratic-final}
\end{align}
Since
\[
(1+|\ell|)^{\tau/2}
\le
(1+|k|+|\ell|)^\tau,
\]
\eqref{eq:quadratic-final} gives \eqref{eq:vector-weyl-main}. The case $R=1$ is absorbed by enlarging $C_{\gamma,\tau}$.
\end{proof}

The same Fourier-mode analysis also records the basic ergodic structure
of the base transformation. We include it here both for completeness
and to emphasize that the quantitative averaging used below does not
come from weak mixing: the skew shift is uniquely ergodic but still has
nontrivial Koopman eigenfunctions.

\begin{lemma}\label{lem:ergodic}
If $\ell\cdot\omega\notin\Z$ for every $0\ne \ell\in\Z^d$, then
$T_\omega$ is uniquely ergodic, with invariant
probability measure equal to Haar measure. It is not weakly mixing.
\end{lemma}
\begin{proof}
The affine map has linear part and inverse
\[
 \begin{pmatrix}I_d&0\\I_d&I_d\end{pmatrix}\in\mathrm{SL}(2d,\Z),
 \qquad T^{-1}(x,y)=(x-\omega,y-x+\omega).
\]
It therefore preserves Haar measure $dz$.

Let $\chi_{k,\ell}(x,y)=e(k\cdot x+\ell\cdot y)$ and
$\mathcal S_Rf=R^{-1}\sum_{j=0}^{R-1}f\circ T^j$.
For $\ell=0$, $k\ne0$, geometric summation gives
\[
 \|\mathcal S_R\chi_{k,0}\|_\infty
 \le\frac{C}{R\|k\cdot\omega\|_\T}\longrightarrow0.
\]
For $\ell\ne0$, \eqref{eq:vdC-normalized}, with $H$ fixed, gives
\begin{align*}
 \|\mathcal S_R\chi_{k,\ell}\|_\infty^2
 &\le\frac2H+\frac4{HR}\sum_{h=1}^{H-1}
       \min\left(R,\frac{C}{\|h\ell\cdot\omega\|_\T}\right),\\
 \limsup_{R\to\infty}\|\mathcal S_R\chi_{k,\ell}\|_\infty^2
 &\le2/H.
\end{align*}
Every denominator is nonzero. Letting $H\to\infty$ proves convergence
to zero for all nonconstant characters. Trigonometric approximation
\cite{Montgomery} now gives, for $f\in C(\T^{2d})$ and any trigonometric
polynomial $p$,
\[
 \|\mathcal S_Rf-\langle f\rangle\|_\infty
 \le2\|f-p\|_\infty+\|\mathcal S_Rp-\langle p\rangle\|_\infty
 \longrightarrow0
\]
after first taking $R\to\infty$ and then $\|f-p\|_\infty\to0$.
If $\nu$ is an invariant probability measure, then
\[
 \int f\,d\nu=\int\mathcal S_Rf\,d\nu\longrightarrow\int f\,dz.
\]
Thus $\nu=dz$, proving unique ergodicity; see also \cite{DF}.
Finally,
\[
 \chi_{k,0}\circ T=e(k\cdot\omega)\chi_{k,0},\qquad k\ne0,
\]
is a nonconstant Koopman eigenfunction, so $T$ is not weakly mixing.
\end{proof}

\subsection{The short-block large-deviation estimate}
\label{subsec:short-block-ldt}

We now combine the vector Weyl bound with the strip-dependent
subharmonic estimates from Sections~\ref{sec:subharmonic} and
\ref{sec:boosting}. The main issue is the anisotropic complex geometry:
the admissible strip has width of order $n^{-1}$ in the $x$-variables
but order one in the $y$-variables. We first regularize in all real
coordinates to obtain a uniform one-sided bound, and then recover a
two-sided deviation estimate by successively averaging in $y$ and
taking $y$-sections at real $x$. The following proposition is the
quantitative short-block input used in the multiscale argument.

\begin{proposition}\label{prop:short-block-LDT}
Fix $\rho_x,\rho_y>0$. Suppose $w_n$ is continuous and separately
subharmonic on $\T_{\rho_x/n}^d\times\T_{\rho_y}^d$, and $|w_n|\le1$.
There exist $A\ge4$, $a>0$, and $C,c,R_0>0$, independent of $n,w_n$,
such that, for $R\ge R_0 n^A$,
\begin{equation}\label{eq:short-block-LDT}
 \mes\left\{z:\abs{\frac1R\sum_{j=0}^{R-1}w_n(T^jz)-\av{w_n}}
 >CR^{-a}\right\}\le Ce^{-cR^a}.
\end{equation}
The constants depend only on $d,\rho_x,\rho_y,\gamma,\tau$, and are
uniform in $\omega\in\DC_d(\gamma,\tau)$, $n$ and $w_n$.
\end{proposition}

\begin{proof}
Set
\[
D_0:=2d,
\qquad
\widetilde\rho_x:=\min\{\rho_x,1/20\},
\qquad
\widetilde\rho_y:=\min\{\rho_y,1/20\}.
\]
We work on
\[
\T_{\widetilde\rho_x/n}^{d}\times
\T_{\widetilde\rho_y}^{d}.
\]
Let $P_t^{(j)}$ denote the periodic Poisson semigroup in the $j$th real
coordinate and put
\[
Q_t:=\prod_{j=1}^{D_0}P_t^{(j)},
\qquad
\widehat{Q_tw_n}(\mathbf m)
=e^{-t|\mathbf m|_1}\widehat w_n(\mathbf m), \quad \mathbf m=(k,l)\in\mathbb{Z}^{2d}.
\]

For the first $d$ coordinates, Lemmas \ref{lem:riesz} and \ref{lem:oneD} give the one-dimensional
cost
\[
K_x\le C\left(\frac{\widetilde\rho_x}{n}\right)^{-4}
\le Cn^4,
\]
whereas for the last $d$ coordinates
\[
K_y\le C\widetilde\rho_y^{-4}\le C.
\]
Positive convolution preserves the strip and the bound.  Iterating the
one-sided estimate coordinatewise therefore gives, for $0<t<1$,
\begin{equation}\label{eq:Qt-one-sided}
w_n\le Q_tw_n+Cn^4t.
\end{equation}

Write
\[
W_R(z):=\frac1R\sum_{j=0}^{R-1}w_n(T_\omega^jz),
\qquad
\mu_n:=\langle w_n\rangle.
\]
Since $Q_t$ preserves the mean,
\begin{equation}
\frac1R\sum_{j=0}^{R-1}(Q_tw_n)(T_\omega^jz)-\mu_n =
\sum_{\mathbf m\in\Z^{D_0}\setminus\{0\}}
 e^{-t|\mathbf m|_1}\widehat w_n(\mathbf m)
 \frac1R\sum_{j=0}^{R-1}
 e\!\left(\mathbf m\cdot T_\omega^jz\right).
\label{eq:regularized-fourier}
\end{equation}
Here $|\widehat w_n(\mathbf m)|\le1$.  Lemma~\ref{lem:vector-weyl} yields, with $\delta:=\frac1{2(\tau+1)}$,
\[
\left|
\frac1R\sum_{j=0}^{R-1}(Q_tw_n)(T_\omega^jz)-\mu_n
\right|
\le
C_{\gamma,\tau}R^{-\delta}
\sum_{\mathbf m\ne0}
 e^{-t|\mathbf m|_1}(1+|\mathbf m|_1)^\tau.
\]
For
\[
\mathcal S_q:=\{\mathbf m\in\Z^{D_0}:|\mathbf m|_1=q\},
\qquad
|\mathcal S_q|\le C_{D_0}(1+q)^{D_0-1},
\]
we obtain
\[
\sum_{\mathbf m\ne0}e^{-t|\mathbf m|_1}(1+|\mathbf m|_1)^\tau
\le
C\sum_{q=1}^{\infty}e^{-tq}(1+q)^{D_0-1+\tau}\le
Ct^{-(D_0+\tau)}, \qquad 0<t<1.
\]
Hence
\[
\sup_z
\left|
\frac1R\sum_{j=0}^{R-1}(Q_tw_n)(T_\omega^jz)-\mu_n
\right|
\le
CR^{-\delta}t^{-(D_0+\tau)}.
\]
Averaging \eqref{eq:Qt-one-sided} along the orbit gives
\[
W_R(z)-\mu_n
\le
C\left(
 n^4t+R^{-\delta}t^{-(D_0+\tau)}
\right),
\qquad z\in\T^{2d}.
\]

Set
\begin{equation}\label{eq:parameters-43}
b_0:=\frac{\delta}{2(D_0+\tau+1)},
\qquad
s_0:=\frac{b_0}{4},
\qquad
\kappa:=\frac{s_0}{8^{D_0+2}},
\end{equation}
and choose
\[
A\ge
\max\left\{4,\frac8\kappa\right\}.
\]
Then
\[
-\delta+b_0(D_0+\tau)<-\frac\delta2,
\qquad
\frac4A\le\frac\kappa2<\frac{b_0}{2}=2s_0.
\]
For $R\ge R_0n^A$ and $t=R^{-b_0}$,
\begin{align}
n^4t
&\le
R^{-b_0+4/A}
\le
R^{-b_0/2}
=R^{-2s_0},
\label{eq:first-error-43}\\
R^{-\delta}t^{-(D_0+\tau)}
&=
R^{-\delta+b_0(D_0+\tau)}
\le
R^{-\delta/2}
\le
R^{-2s_0}.
\label{eq:second-error-43}
\end{align}
Increasing $R_0$ once,
\begin{equation}\label{eq:uniform-upper-43}
W_R(x,y)-\mu_n\le R^{-s_0},
\qquad (x,y)\in\T^{2d}.
\end{equation}

The one-dimensional strip costs satisfy
\[
K_x\le Cn^4
\le CR^{4/A}
\le CR^{\kappa/2},
\qquad
K_y\le C\le CR^\kappa.
\]
By Lemma~\ref{lem:positive-average} the same bounds hold after any finite sequence of real
coordinate averages.  Thus all sections used below have analytic cost at
most $CR^\kappa$.

Define
\[
\overline w_n(x):=\int_{\T^d}w_n(x,y)\,dy.
\]
For real $x$, invariance of Haar measure under real translations gives
\[
G_R(x)
:=\int_{\T^d}\bigl(W_R(x,y)-\mu_n\bigr)\,dy
=\frac1R\sum_{j=0}^{R-1}
\overline w_n(x+j\omega)-\mu_n.
\]
The right-hand side defines the complex extension of $G_R$ to
$\T_{\widetilde\rho_x/n}^{d}$.  By Lemma~\ref{lem:positive-average},
\begin{equation}\label{eq:GR-properties}
|G_R|\le2,
\qquad
G_R\ \text{is separately subharmonic},
\qquad
\langle G_R\rangle_x=0.
\end{equation}
Integrating \eqref{eq:uniform-upper-43} in $y$ gives $G_R(x)\le R^{-s_0}$, with $x\in\T^d.$ Moreover,
\[
\kappa
=\frac{s_0}{8^{2d+2}}
\le
\frac{s_0}{8^{d+1}}.
\]
Lemma~\ref{lem:one-sided-concentration}, with scale parameter $R$, dimension $d$, exponent $s_0$,
and cost $CR^\kappa$, yields
\begin{equation}\label{eq:GR-LDT}
\mes_x\left\{
|G_R(x)|>CR^{-s_1}
\right\}
\le
C\exp(-cR^{s_1}),
\qquad
s_1:=\frac{s_0}{8^d}.
\end{equation}
Set
\[
E_R:=\{x\in\T^d:|G_R(x)|>CR^{-s_1}\}.
\]

Fix $x\in\T^d\setminus E_R$ and set
\[
F_x(y):=W_R(x,y)-\mu_n.
\]
Since $x$ is real,
\[
T_\omega^j(x,y)
=
\left(
 x+j\omega,
 y+jx+\frac{j(j-1)}2\omega
\right),
\]
and the $y$-variable is translated only by real vectors.  Hence
$F_x$ extends separately subharmonically to
$\T_{\widetilde\rho_y}^{d}$ and
\[
|F_x|\le2,
\qquad
F_x(y)\le R^{-s_0}\le R^{-s_1},
\qquad
|\langle F_x\rangle_y|=|G_R(x)|\le CR^{-s_1}.
\]
The $y$-section costs are bounded by $CR^\kappa$, and
\[
\kappa
=\frac{s_0}{8^{2d+2}}
=\frac{s_1}{8^{d+2}}
\le
\frac{s_1}{8^{d+1}}.
\]
Applying Lemma~\ref{lem:one-sided-concentration} again, now in dimension $d$ with exponent $s_1$, gives
\begin{equation}\label{eq:y-fiber-LDT}
\mes_y\left\{
|F_x(y)|>CR^{-a}
\right\}
\le
C\exp(-cR^a),
\qquad
 a:=\frac{s_1}{8^d}=\frac{s_0}{8^{2d}}>0,
\end{equation}
uniformly for $x\notin E_R$.
Therefore, by \eqref{eq:GR-LDT}, \eqref{eq:y-fiber-LDT}, and Fubini's theorem,
\begin{align}
\mes\left\{
|W_R-\mu_n|>CR^{-a}
\right\}
&\le
\mes_x(E_R)
+
\int_{\T^d\setminus E_R}
\mes_y\{|F_x|>CR^{-a}\}\,dx
\notag\\
&\le
C e^{-cR^{s_1}}+C e^{-cR^a}
\le
C e^{-cR^a}.
\end{align}
This is \eqref{eq:short-block-LDT}.
\end{proof}

\begin{remark}[Complex domains and the short-block regime]
The full averaged function \(W_R\) is not claimed to admit a uniformly bounded extension to the original product strip. Its $y$-fibers and its $y$-average have the separate
extensions just described. These are all the extensions the proof uses.
For fixed $n$ the result implies every prescribed polynomial exceptional
measure bound, but uniformity in $n$ requires $R\ge R_0n^A$.
\end{remark}

% Section 5
\section{Transfer matrices and the large-coupling initial scale}
\label{sec:initial}

\subsection{Analytic normalization and uniform sublevel estimates}
\label{subsec:initial-preparation}

The initial-scale argument requires two analytic inputs that are uniform
in the real energy.  The first places the normalized logarithmic transfer
matrix norm on the anisotropic complex strip dictated by the skew-shift
shear and records the plurisubharmonic structure needed for
Proposition~\ref{prop:short-block-LDT}.  The second is a uniform
transversality estimate for the level sets of the analytic sampling
function.  We establish both ingredients before turning to the
large-coupling initial scale.

\begin{lemma}\label{lem:analytic-normalization}
For the choice \eqref{eq:S}, increasing the fixed constant $C_V$ if necessary,
\begin{equation}\label{eq:analytic-normalization}
0\le u_n(z;\lambda,E)\le S(\lambda,E)
\qquad
\text{on }
\T_{\rho/(4n)}^d\times\T_{\rho/4}^d .
\end{equation}
Moreover, $u_n$ is continuous and plurisubharmonic on this domain;
in particular, it is separately subharmonic.  Consequently,
$u_n/S$ satisfies the hypotheses of Proposition~\ref{prop:short-block-LDT} with
$\rho_x=\rho_y=\rho/4$, uniformly in real $E$ and $\lambda>0$.
\end{lemma}

\begin{proof}
Put
\[
\Omega_n
:=
\T_{\rho/(4n)}^d\times\T_{\rho/4}^d,
\qquad
M_\rho
:=
\sup_{z\in\T_{\rho/2}^{2d}}|V(z)|<\infty.
\]
Choose $C_V$ so that
\begin{equation}\label{eq:CV-choice}
C_V\ge 2(1+M_\rho),
\qquad
\log(2C_V)\ge1.
\end{equation}

Write
\[
z=(x,y),
\qquad
x=x_0+i\xi,
\qquad
y=y_0+i\eta,
\qquad z\in\Omega_n.
\]
For $0\le j<n$, since $\omega\in\T^d$ is real,
\begin{equation}
\norm{\Im(x+j\omega)}_\infty=
\norm{\xi}_\infty
<\frac{\rho}{4n},
\label{eq:x-strip}
\end{equation}
and
\begin{equation}
\norm{\Im\!\left(y+jx+\frac{j(j-1)}2\omega\right)}_\infty=
\norm{\eta+j\xi}_\infty\le
\norm{\eta}_\infty+j\norm{\xi}_\infty
<
\frac\rho4+\frac{j\rho}{4n}
<
\frac\rho2.
\label{eq:y-strip}
\end{equation}
Hence
\begin{equation}\label{eq:orbit-inside}
T_\omega^j(\Omega_n)
\subset
\T_{\rho/2}^{2d}
\subset
\T_\rho^{2d},
\qquad
0\le j<n.
\end{equation}

For $z\in\Omega_n$ and $0\le j<n$, using
\eqref{eq:orbit-inside} and \eqref{eq:CV-choice},
\begin{align}
\norm{A(T_\omega^jz;\lambda,E)}
&\le
\norm{A(T_\omega^jz;\lambda,E)}_{\mathrm F}
\notag\\
&=
\left(
|E-\lambda V(T_\omega^jz)|^2+2
\right)^{1/2}
\notag\\
&\le
|E|+\lambda M_\rho+2
\notag\\
&\le
C_V(2+\lambda+|E|)
=e^{S(\lambda,E)}.
\label{eq:one-step-bound}
\end{align}
Here \(\|\cdot\|_{\mathrm F}\) denotes the Frobenius norm. Therefore
\begin{equation}\label{eq:Mn-upper}
\norm{M_n(z;\lambda,E)}
\le
\prod_{j=0}^{n-1}
\norm{A(T_\omega^jz;\lambda,E)}
\le
e^{nS(\lambda,E)}.
\end{equation}
Moreover,
\begin{equation}\label{eq:det-one}
\det A(z;\lambda,E)=1,
\qquad
\det M_n(z;\lambda,E)=1.
\end{equation}
Let
\[
s_1(z)\ge s_2(z)>0
\]
be the singular values of $M_n(z;\lambda,E)$.  Then
\begin{equation}\label{eq:singular-values}
s_1(z)s_2(z)
=
|\det M_n(z;\lambda,E)|
=1,
\qquad
s_1(z)=\norm{M_n(z;\lambda,E)}\ge1.
\end{equation}
Combining \eqref{eq:Mn-upper} and \eqref{eq:singular-values},
\begin{equation}\label{eq:un-bounds}
1
\le
\norm{M_n(z;\lambda,E)}
\le
e^{nS(\lambda,E)},
\qquad
0\le u_n(z;\lambda,E)\le S(\lambda,E).
\end{equation}

By \eqref{eq:orbit-inside}, every entry of $M_n$ is holomorphic on
$\Omega_n$; by \eqref{eq:det-one}, $M_n(z)$ is invertible there.
Hence
\begin{equation}\label{eq:continuity}
z\longmapsto\log\norm{M_n(z;\lambda,E)}
\end{equation}
is continuous on $\Omega_n$.

Fix $z_0\in\Omega_n$, $v\in\C^{2d}$, and $r>0$ such that
\[
z_0+\zeta v\in\Omega_n,
\qquad |\zeta|\le r.
\]
Choose $\xi_0\in\C^2$, $\norm{\xi_0}_2=1$, with
\begin{equation}\label{eq:max-vector}
\norm{M_n(z_0;\lambda,E)\xi_0}_2
=
\norm{M_n(z_0;\lambda,E)}.
\end{equation}
The vector-valued map
\[
F(\zeta)
:=
M_n(z_0+\zeta v;\lambda,E)\xi_0
\]
is holomorphic and nonvanishing.  Thus
\[
\log\norm{F(\zeta)}_2
=
\frac12\log\!
\left(
|F_1(\zeta)|^2+|F_2(\zeta)|^2
\right)
\]
is subharmonic.  By the submean inequality and
\eqref{eq:max-vector},
\begin{align}
\log\norm{M_n(z_0;\lambda,E)}
&=
\log\norm{F(0)}_2
\notag\\
&\le
\frac1{2\pi}
\int_0^{2\pi}
\log\norm{F(re^{i\theta})}_2\,d\theta
\notag\\
&\le
\frac1{2\pi}
\int_0^{2\pi}
\log\norm{M_n(z_0+re^{i\theta}v;\lambda,E)}\,d\theta.
\label{eq:line-submean}
\end{align}
Thus \eqref{eq:continuity} is subharmonic on every complex line in
$\Omega_n$ and hence plurisubharmonic.  Division by the positive
constant $n$ gives the same conclusion for $u_n$.

Finally, by \eqref{eq:CV-choice}, $S\ge1$, and \eqref{eq:un-bounds}
gives
\[
0\le \frac{u_n}{S}\le1
\qquad\text{on}\qquad
\T_{(\rho/4)/n}^{d}\times\T_{\rho/4}^{d}.
\]
Therefore $u_n/S$ is continuous and separately subharmonic on the
precise anisotropic strip required in Proposition~\ref{prop:short-block-LDT}, with constants
independent of $E\in\R$ and $\lambda>0$.
\end{proof}

The preceding lemma supplies the complex-analytic control needed for
short-block averaging.  To obtain a large-coupling lower bound at the
initial scale, we also need to exclude phases for which one of the
diagonal factors is anomalously small.  This is provided by the following
uniform sublevel estimate, a compact analytic transversality argument
related to \cite[Lemma~11.4]{GS} and the smooth multivariable estimates
of \cite{Klein}.  We include the proof in order to make explicit that the
constants are uniform in the translating level.

\begin{lemma}\label{lem:uniform-sublevel}
Let $V:\T^{2d}\to\R$ be nonconstant and real analytic. There exist
constants $C_{\rm sub}>0, \beta\in(0,1],$ depending only on $V$, such that, for every $t\in\R$ and every
$0<\varepsilon<1$,
\begin{equation}\label{eq:uniform-sublevel}
 \mes\bigl\{z\in\T^{2d}: |V(z)-t|<\varepsilon\bigr\}
 \le C_{\rm sub}\varepsilon^{\beta}.
\end{equation}
\end{lemma}

\begin{proof}
Let $I\subset\R$ be a bounded interval, $k\ge1$, and
$g\in C^k(I)$ satisfy
\begin{equation}\label{eq:derivative-lower}
 \inf_{s\in I}|g^{(k)}(s)|\ge c>0.
\end{equation}
For $t\in\R$ and $\varepsilon>0$, set
\[
 E:=\{s\in I: |g(s)-t|<\varepsilon\},
 \qquad
 a:=|E|.
\]
If $a=0$, there is nothing to prove. Assume $a>0$ and put
\begin{equation}\label{eq:delta-selection}
 \delta:=\frac{a}{2(k+1)}.
\end{equation}
There exist $k+1$ points
\[
 s_0<\cdots<s_k,
 \qquad
 s_j\in E,
 \qquad
 |s_i-s_j|\ge\delta\quad(i\ne j).
\]
Indeed, a maximal $\delta$-separated subset of $E$ with at most $k$
points would cover $E$ by at most $k$ intervals of length $2\delta$,
and hence
\[
 a=|E|\le 2k\delta
   =\frac{k}{k+1}a<a,
\]
a contradiction.

The divided-difference mean-value theorem (see \cite{deBoor}) gives some
$\xi\in(s_0,s_k)$ such that
\[
 \frac{g^{(k)}(\xi)}{k!}
 =\sum_{i=0}^k
 \frac{g(s_i)-t}{\prod_{j\ne i}(s_i-s_j)}.
\]
Using \eqref{eq:derivative-lower}, $|g(s_i)-t|<\varepsilon$, and
\eqref{eq:delta-selection},
\begin{equation}
 \frac{c}{k!}\le
 \sum_{i=0}^k
 \frac{\varepsilon}
 {\prod_{j\ne i}|s_i-s_j|}\le
 (k+1)\varepsilon\delta^{-k}
 =
 (k+1)\varepsilon
 \left(\frac{2(k+1)}{a}\right)^k.
 \label{eq:sublevel-one-dimensional-computation}
\end{equation}
Therefore
\begin{equation}\label{eq:sublevel-one-dimensional}
 |E|
 \le
 C_k\left(\frac{\varepsilon}{c}\right)^{1/k},
 \qquad
 C_k:=2(k+1)\bigl(k!(k+1)\bigr)^{1/k}.
\end{equation}
The constant in \eqref{eq:sublevel-one-dimensional} is independent of
the translating level $t$.

Fix $z_0\in\T^{2d}$ and use a real-analytic lift near $z_0$. For
$k\ge1$ define the homogeneous polynomial
\[
 P_{k,z_0}(\xi)
 :=
 \left.\frac{d^k}{ds^k}V(z_0+s\xi)\right|_{s=0}
 =\partial_\xi^kV(z_0),
 \qquad
 \xi\in\R^{2d}.
\]
If $ P_{k,z_0}\equiv0$ for $k\ge1$, then every positive-order Taylor coefficient of $V$ at $z_0$ vanishes.
Analyticity gives $ V\equiv V(z_0)$ on a neighborhood of $z_0$, hence on the connected torus by analytic
continuation, contradicting the hypothesis. Thus, for every $z_0$,
there exist
\[
 k(z_0)\ge1,
 \qquad
 \xi(z_0)\in\R^{2d},
 \qquad
 |\xi(z_0)|_2=1,
\]
such that $\partial_{\xi(z_0)}^{k(z_0)}V(z_0)\ne0.$ By continuity, after shrinking a local Euclidean box $B_{z_0}$,
\[
 \bigl|\partial_{\xi(z_0)}^{k(z_0)}V(z)\bigr|
 \ge c(z_0)>0,
 \qquad
 z\in B_{z_0}.
\]
Compactness of $\T^{2d}$ yields finitely many boxes
\[
 \T^{2d}\subset\bigcup_{i=1}^J B_i,
\]
with associated data
\begin{equation}\label{eq:finite-data}
 k_i\ge1,
 \qquad
 |\xi_i|_2=1,
 \qquad
 c_i>0,
 \qquad
 |\partial_{\xi_i}^{k_i}V|\ge c_i
 \quad\text{on }B_i.
\end{equation}
Choose each $B_i$ in orthogonal coordinates with first axis $\xi_i$:
\[
 B_i=
 \{z_i+s\xi_i+y:
   s\in I_i,\ y\in Q_i\subset\xi_i^\perp\},
\]
where $I_i$ is an interval and $Q_i$ is a bounded rectangle.

For fixed $y\in Q_i$, define
\[
 g_{i,y}(s):=V(z_i+s\xi_i+y),
 \qquad
 s\in I_i.
\]
By \eqref{eq:finite-data},
\[
 |g_{i,y}^{(k_i)}(s)|
 =
 |\partial_{\xi_i}^{k_i}V(z_i+s\xi_i+y)|
 \ge c_i.
\]
Hence \eqref{eq:sublevel-one-dimensional} gives, uniformly in
$y\in Q_i$ and $t\in\R$,
\[
 \left|
 \{s\in I_i:
   |V(z_i+s\xi_i+y)-t|<\varepsilon\}
 \right|
 \le
 C_{k_i}c_i^{-1/k_i}\varepsilon^{1/k_i}.
\]
Fubini's theorem therefore yields
\begin{equation}\label{eq:box-sublevel}
 \mes\bigl(B_i\cap\{|V-t|<\varepsilon\}\bigr)
 \le
 |Q_i|C_{k_i}c_i^{-1/k_i}\varepsilon^{1/k_i}.
\end{equation}
Set
\begin{equation}\label{eq:global-beta}
 k_*:=\max_{1\le i\le J}k_i,
 \qquad
 \beta:=\frac1{k_*}.
\end{equation}
Since $0<\varepsilon<1$ and $k_i\le k_*$,
\[
 \varepsilon^{1/k_i}
 \le
 \varepsilon^{1/k_*}
 =\varepsilon^\beta.
\]
Summing \eqref{eq:box-sublevel} over the finite covering gives
\[
 \mes\{|V-t|<\varepsilon\}\le
 \sum_{i=1}^J
 |Q_i|C_{k_i}c_i^{-1/k_i}
 \varepsilon^{1/k_i}\le
 C_{\rm sub}\varepsilon^\beta,
\]
where
\[
 C_{\rm sub}
 :=
 \sum_{i=1}^J
 |Q_i|C_{k_i}c_i^{-1/k_i}.
\]
All quantities in \eqref{eq:finite-data}--\eqref{eq:global-beta}
depend only on $V$. In particular, $C_{\rm sub}$ and $\beta$ are
independent of $t$ and $\varepsilon$.
\end{proof}

\subsection{The large-coupling initial scale}
\label{subsec:large-coupling-initial}

We now combine the anisotropic normalization with the uniform sublevel
bound.  For energies comparable with the coupling, the bad set is formed
by phases for which some diagonal entry approaches zero; outside this set,
a finite-volume determinant factorization gives a pointwise lower bound
for the transfer matrix.  For energies much larger than the coupling,
direct diagonal dominance gives the same conclusions without an
exceptional set.  The two regimes yield a uniform initial-scale LDT,
positivity, and the first two-scale gap estimate.

\begin{lemma}\label{lem:initial-scale}
Fix $\epsilon>0$ small and $Q>0$. There exist
\[
B=B(\eps,Q,V)>0,
\qquad
\lambda_*=\lambda_*(\eps,Q,V)>1,
\]
such that for every integer $n_0\ge2$, every $\lambda\ge \lambda_* n_0^B,$ and every $E\in\mathbb R$,
\begin{align}
&\mes\{|u_{n_0}-L_{n_0}|>2\eps S\}
 +\mes\{|u_{2n_0}-L_{2n_0}|>2\eps S\}
 \le n_0^{-Q}, \label{eq:initial-ldt}\\
&\min\{L_{n_0},L_{2n_0}\}\ge(1-2\eps)S, \label{eq:initial-positive}\\
&0\le L_{n_0}-L_{2n_0}\le2\eps S. \label{eq:initial-gap}
\end{align}
The constants are uniform in the real energy and in the frequency in the fixed class.
\end{lemma}

\begin{proof}
Put $V_*:=\norm{V}_{L^\infty(\T^{2d})}.$ Choose
\begin{equation}\label{eq:zeta-choice}
0<\zeta\le \min\left\{\frac18,\frac{\eps}{32}\right\},
\end{equation}
and then choose $B$ so that
\begin{equation}\label{eq:B-choice}
B\beta\zeta>Q+2.
\end{equation}
Only $\lambda_*$ will be enlarged below.

We first treat the energy range in which the diagonal terms can become
small on a nonempty exceptional set.  The sublevel estimate controls the
measure of this set uniformly in the energy.

Assume
\begin{equation}\label{eq:moderate-energy}
|E|\le2\lambda(V_*+1).
\end{equation}
For $m\in\{n_0,2n_0\}$ define
\[
\mathcal B_m(E)
:=\bigcup_{j=0}^{m-1}
\left\{z:\left|V(T_\omega^jz)-\frac{E}{\lambda}\right|<\lambda^{-\zeta}\right\}.
\]
Haar invariance and Lemma~\ref{lem:uniform-sublevel} give
\begin{align}
\mes\mathcal B_m(E)
&\le \sum_{j=0}^{m-1}
\mes\left\{z:\left|V(T_\omega^jz)-\frac{E}{\lambda}\right|<\lambda^{-\zeta}\right\}\notag\\
&=m\,\mes\left\{z:\left|V(z)-\frac{E}{\lambda}\right|<\lambda^{-\zeta}\right\}\notag\\
&\le \Csub m\lambda^{-\beta\zeta}. \label{eq:bad-measure}
\end{align}
Hence, for $m\in\{n_0,2n_0\}$,
\[
\mes\mathcal B_{n_0}(E)+\mes\mathcal B_{2n_0}(E)
\le C n_0\lambda^{-\beta\zeta}.
\]

Let $H_{[0,m-1]}(z)$ be the Dirichlet restriction of the Schr\"odinger operator to $[0,m-1]$, and set
\[
D_m(z):=\operatorname{diag}\bigl(\lambda V(T_\omega^jz)-E\bigr)_{j=0}^{m-1}.
\]
Let $J_m$ be the $m\times m$ matrix with ones on the two adjacent diagonals and zeros elsewhere. Then
\begin{equation}\label{eq:finite-volume-matrix}
D_m+J_m=-(E-H_{[0,m-1]}),
\qquad
\norm{J_m}\le2.
\end{equation}
For $z\notin\mathcal B_m(E)$,
\begin{equation}\label{eq:diagonal-lower}
|\lambda V(T_\omega^jz)-E|
\ge \lambda^{1-\zeta},
\qquad 0\le j<m,
\end{equation}
and therefore
\[
\norm{D_m^{-1}J_m}
\le2\lambda^{-(1-\zeta)}
=:\eta_\lambda.
\]
Increase $\lambda_*$ so that $\eta_\lambda\le1/2$.

Define
\[
p_m(z,E):=\det(E-H_{[0,m-1]}(z)),
\qquad p_{-1}=0,\quad p_0=1.
\]
Expansion along the last row gives
\[
p_m
=\bigl(E-\lambda V(T_\omega^{m-1}z)\bigr)p_{m-1}-p_{m-2}.
\]
The transfer recursion gives the same recurrence and initial data, hence
\begin{equation}\label{eq:M11}
(M_m(z;\lambda,E))_{11}=p_m(z,E).
\end{equation}
By \eqref{eq:finite-volume-matrix},
\[
|p_m|=|\det(D_m+J_m)|=|\det D_m|\,|\det(I+D_m^{-1}J_m)|.
\]
Since
\[
s_{\min}(I+D_m^{-1}J_m)
\ge1-\norm{D_m^{-1}J_m}
\ge1-\eta_\lambda,
\]
we obtain
\[
|p_m|\ge
\prod_{j=0}^{m-1}|\lambda V(T_\omega^jz)-E|
\,(1-\eta_\lambda)^m\ge
\lambda^{m(1-\zeta)}(1-\eta_\lambda)^m.
\]
Using \eqref{eq:M11}, for $z\notin\mathcal B_m(E)$, we have
\begin{equation}\label{eq:u-lower}
u_m(z)
=\frac1m\log\norm{M_m(z)}
\ge\frac1m\log|p_m|
\ge(1-\zeta)\log\lambda+\log(1-\eta_\lambda).
\end{equation}

By Lemma~\ref{lem:analytic-normalization},
\begin{equation}\label{eq:u-upper}
0\le u_m(z)\le S.
\end{equation}
Under \eqref{eq:moderate-energy},
\[
S=\log\!\bigl(C_V(2+\lambda+|E|)\bigr)\le \log\lambda+C_0(V),
\]
while $S\ge\log\lambda$. Hence, off $\mathcal B_m(E)$,
\begin{align}
0\le1-\frac{u_m}{S}
&\le
\frac{C_0(V)+\zeta\log\lambda-\log(1-\eta_\lambda)}{S}. \label{eq:relative-defect}
\end{align}
By \eqref{eq:zeta-choice}, after increasing $\lambda_*$,
\begin{equation}\label{eq:pointwise-relative}
0\le1-\frac{u_m}{S}\le\frac{\eps}{4},
\qquad z\notin\mathcal B_m(E),
\quad m\in\{n_0,2n_0\}.
\end{equation}

From \eqref{eq:B-choice} and $\lambda\ge\lambda_*n_0^B,$ we have
\[
Cn_0\lambda^{-\beta\zeta}\le C\lambda_*^{-\beta\zeta}
   n_0^{1-B\beta\zeta}\le C\lambda_*^{-\beta\zeta}n_0^{-Q-1}.
\]
Increase $\lambda_*$ so that
\begin{equation}\label{eq:bad-small}
\mes\mathcal B_{n_0}(E)+\mes\mathcal B_{2n_0}(E)
\le \min\left\{\frac{\eps}{8},n_0^{-Q}\right\}.
\end{equation}
Using \eqref{eq:u-upper}, \eqref{eq:pointwise-relative}, and $u_m\ge0$,
\begin{equation}\label{eq:L-relative}
0\le1-\frac{L_m}{S}=\int_{\T^{2d}}\left(1-\frac{u_m}{S}\right)dz\le \frac{\eps}{4}+\mes\mathcal B_m(E)
\le\frac{3\eps}{8}<\frac{\eps}{2}.
\end{equation}
Consequently, for $z\notin\mathcal B_m(E)$,
\[
|u_m(z)-L_m|\le (S-u_m(z))+(S-L_m)
\le \frac{3\eps}{4}S
<2\eps S.
\]
Thus
\[
\{|u_m-L_m|>2\eps S\}\subset\mathcal B_m(E),
\qquad m\in\{n_0,2n_0\},
\]
and \eqref{eq:initial-ldt} follows from \eqref{eq:bad-small}.

Moreover, \eqref{eq:L-relative} gives
\[
\min\{L_{n_0},L_{2n_0}\}
\ge\left(1-\frac{\eps}{2}\right)S
\ge(1-2\eps)S.
\]
Subadditivity yields
\begin{equation}\label{eq:subadditivity-gap}
0\le L_{n_0}-L_{2n_0},
\end{equation}
and therefore
\begin{equation}\label{eq:gap-upper}
0\le L_{n_0}-L_{2n_0}
\le S-L_{2n_0}
\le\frac{\eps}{2}S
\le2\eps S.
\end{equation}
This proves \eqref{eq:initial-positive}--\eqref{eq:initial-gap} in the range \eqref{eq:moderate-energy}.

It remains to consider energies outside the coupling scale.  In this
regime every diagonal entry is uniformly large, so the same determinant
factorization becomes global and no exceptional set is needed.

Assume now
\begin{equation}\label{eq:large-energy}
|E|>2\lambda(V_*+1).
\end{equation}
Then, for every real phase and every $j$,
\[
\frac{|E|}{2}
\le |\lambda V(T_\omega^jz)-E|
\le\frac{3|E|}{2}.
\]
Hence
\[
\norm{D_m^{-1}J_m}\le\frac4{|E|}.
\]
After increasing $\lambda_*$, $4/|E|\le1/2$. The same determinant factorization gives
\begin{align}
|p_m|
&\ge
\left(\frac{|E|}{2}\right)^m
\left(1-\frac4{|E|}\right)^m,\notag\\
u_m(z)
&\ge
\log|E|-\log2+\log\left(1-\frac4{|E|}\right). \label{eq:large-u-lower}
\end{align}
Also
\[
\norm{A(T_\omega^jz;\lambda,E)}
\le |E-\lambda V(T_\omega^jz)|+2
\le 2|E|,
\]
for $\lambda_*$ sufficiently large. Thus
\[
|u_m(z)-\log|E||\le C_1(V),
\qquad m\in\{n_0,2n_0\}.
\]
Averaging,
\[
|L_m-\log|E||\le C_1(V).
\]
Moreover, \eqref{eq:large-energy} implies
\[
\lambda\le \frac{|E|}{2(V_*+1)},
\]
and hence
\[
0\le S-\log|E|\le C_2(V).
\]
Therefore, with $C_3=C_1+C_2$,
\[
|u_m-S|+|L_m-S|\le2C_3.
\]
Since $S\ge\log\lambda$ and $\lambda\ge\lambda_*n_0^B$, enlarge $\lambda_*$ so that
\begin{equation}\label{eq:large-absorb}
2C_3\le\frac{\eps}{2}S.
\end{equation}
Then, for both $m=n_0,2n_0$,
\begin{align}
|u_m-L_m|&\le |u_m-S|+|L_m-S|\le\frac{\eps}{2}S<2\eps S,\notag\\
L_m&\ge(1-\eps/2)S. \label{eq:large-conclusions}
\end{align}
The exceptional sets are empty. Finally,
\begin{equation}\label{eq:large-gap}
0\le L_{n_0}-L_{2n_0}
\le |L_{n_0}-S|+|L_{2n_0}-S|
\le\frac{\eps}{2}S
\le2\eps S.
\end{equation}
This proves all three assertions in the range \eqref{eq:large-energy} and completes the proof.
\end{proof}

\section{An averaged avalanche principle and one induction step}
\label{sec:AP}

\subsection{The all-offset avalanche comparison}
\label{subsec:averaged-avalanche}

The purpose of this subsection is to convert finite-scale control at the
two adjacent scales $n$ and $2n$ into a deterministic comparison at a
larger scale $N$.  The key point is to apply the avalanche principle to
every residue class modulo $n$ and then average over the offsets.  The
resulting boundary error is of order $n/N$, while the exceptional set is
controlled by the finite-scale deviation probability.

We use the following standard $2\times2$ avalanche principle of
Goldstein--Schlag \cite{GS}, also stated in \cite[Proposition 2.1]{BGS}:
if $B_1,\ldots,B_m\in\mathrm{SL}(2,\C)$ satisfy
\begin{equation}\label{eq:APhyp}
 \min_j\norm{B_j}\ge\mu\ge m,
 \qquad
 \max_{j<m}\big(\log\norm{B_{j+1}}+\log\norm{B_j}
 -\log\norm{B_{j+1}B_j}\big)\le\tfrac12\log\mu,
\end{equation}
then
\begin{equation}\label{eq:AP}
 \abs{\log\norm{B_m\cdots B_1}
 +\sum_{j=2}^{m-1}\log\norm{B_j}
 -\sum_{j=1}^{m-1}\log\norm{B_{j+1}B_j}}
 \le Cm/\mu.
\end{equation}
This classical matrix result is an external input; its hypotheses will be
verified quantitatively below.

For this section fix $\lambda,E$ and write
\begin{equation}
 h_j=u_j/S,\quad \ell_j=L_j/S,\quad
 \gamma_*=1/2,\quad \varepsilon_* =\gamma_*/100.
\end{equation}
Then $0\le h_j\le1$ on the real torus.

Before applying the avalanche principle, we record the elementary
counting identities that make the averaging over all block offsets exact.
They ensure that the interior $n$- and $2n$-blocks are counted once after
summing over the residue classes modulo $n$.

\begin{lemma}\label{lem:offset-counting}
Let $N\ge 4n$ be integers and, for $0\le r<n$, set
\[
    m_r:=\left\lfloor \frac{N-r}{n}\right\rfloor .
\]
For arbitrary scalar families $(a_k)$ and $(b_k)$,
\begin{align}
\sum_{r=0}^{n-1}\sum_{j=0}^{m_r-2} a_{r+jn}
    &=\sum_{k=0}^{N-2n} a_k, \label{eq:offset-a}\\
\sum_{r=0}^{n-1}\sum_{j=1}^{m_r-2} b_{r+jn}
    &=\sum_{k=n}^{N-2n} b_k. \label{eq:offset-b}
\end{align}
Empty sums are understood to be zero.
\end{lemma}

\begin{proof}
For every \(k\in\mathbb Z_{\ge0}\), Euclidean division by \(n\) gives the unique representation
\[
k=r+jn,
\qquad
0\le r<n,
\qquad
j\ge0,
\]
with
\[
j=\left\lfloor\frac{k}{n}\right\rfloor,
\qquad
r=k-n\left\lfloor\frac{k}{n}\right\rfloor.
\]
For such a pair, if $ j\le m_r-2$, we have $k\le N-2n$.
Hence, when $(r,j)\mapsto r+jn$,
\[
    \bigl\{(r,j):0\le r<n,\ 0\le j\le m_r-2\bigr\}
    \rightarrow\{0,1,\ldots,N-2n\}
\]
is a bijection. Therefore
\[
    \sum_{r=0}^{n-1}\sum_{j=0}^{m_r-2}a_{r+jn}
    =\sum_{k=0}^{N-2n}a_k,
\]
which proves \eqref{eq:offset-a}.

Moreover, since $0\le r<n$, $j\ge1$, we have $k\ge n$. When $(r,j)\mapsto r+jn$,
\[
    \bigl\{(r,j):0\le r<n,\ 1\le j\le m_r-2\bigr\}
    \rightarrow
    \{n,n+1,\ldots,N-2n\}
\]
is a bijection. Thus
\[
    \sum_{r=0}^{n-1}\sum_{j=1}^{m_r-2}b_{r+jn}
    =\sum_{k=n}^{N-2n}b_k,
\]
and \eqref{eq:offset-b} follows.
\end{proof}

\medskip
\noindent\textit{Remark.}
The counting identities themselves require only $N\ge2n$; the stronger hypothesis
$N\ge4n$ is retained because it is the scale separation used in the subsequent
all-offset avalanche estimate.

We now combine this exact counting with the avalanche principle.  The
finite-scale LDT enters only through the bad set on which one of the
$n$- or $2n$-block estimates fails.  Outside that set, the norm lower
bound and the small two-scale gap verify the avalanche hypotheses
uniformly for every offset.

\begin{lemma}\label{lem:APaverage}
Suppose
\begin{equation}\label{eq:APscalebounds}
 \ell_n,\ell_{2n}\ge\gamma_*,\qquad
 0\le\ell_n-\ell_{2n}\le\gamma_*/40,
\end{equation}
and let
\begin{equation}
 q_n=\mes\{|h_n-\ell_n|>\varepsilon_*\}
 +\mes\{|h_{2n}-\ell_{2n}|>\varepsilon_*\}.
\end{equation}
Assume
\begin{equation}\label{eq:Nn-hyp}
N\ge4n,
\qquad
nS\ge100\log(4N).
\end{equation}
There is a set $\mathcal B$ with $\mes\mathcal B\le CNq_n$ such that, off $\mathcal B$,
\begin{equation}\label{eq:APavg}
 \abs{h_N+\frac1N\sum_{j=0}^{N-1}h_n\circ T^j
 -\frac2N\sum_{j=0}^{N-1}h_{2n}\circ T^j}
 \le Cn/N.
\end{equation}
The analogous assertion at scale $2N$ holds after at most the same order
of enlargement of $\mathcal B$. Consequently,
\begin{align}
 |\ell_N+\ell_n-2\ell_{2n}|&\le C(n/N+Nq_n),\label{eq:APmean}\\
 0\le\ell_N-\ell_{2N}&\le C(n/N+Nq_n).\label{eq:APgap}
\end{align}
\end{lemma}

\begin{proof}
Put
\[
\mathcal E_j:=\{|h_j-\ell_j|>\eps_*\},
\qquad j\in\{n,2n\},
\]
and define
\begin{equation}\label{eq:Bdef}
\mathcal B
:=
\bigcup_{k=0}^{2N}T^{-k}(\mathcal E_n\cup\mathcal E_{2n}).
\end{equation}
Since $T$ preserves Haar measure,
\begin{equation}\label{eq:Bbound}
\mes\mathcal B
\le(2N+1)q_n
\le C Nq_n.
\end{equation}
Fix $z\notin\mathcal B$. For every block start used below,
\[
\log\norm{M_n(T^kz)}
=nS h_n(T^kz)
\ge nS(\ell_n-\eps_*)
\ge nS(\gamma_*-\eps_*).
\]
Set
\begin{equation}\label{eq:mu-def}
\mu:=\exp\bigl((\gamma_*-\eps_*)nS\bigr).
\end{equation}
For two consecutive $n$-blocks,
\begin{align}
&\log\norm{M_n(T^{k+n}z)}
 +\log\norm{M_n(T^kz)}
 -\log\norm{M_{2n}(T^kz)}
\notag\\
&\quad
=nS\bigl(h_n(T^{k+n}z)+h_n(T^kz)-2h_{2n}(T^kz)\bigr)
\notag\\
&\quad
\le
2nS(\ell_n-\ell_{2n})+4nS\eps_*
\notag\\
&\quad
\le
\left(\frac{\gamma_*}{20}+4\eps_*\right)nS
<\frac12(\gamma_*-\eps_*)nS
=\frac12\log\mu.
\label{eq:cancellation}
\end{align}
Moreover, \eqref{eq:Nn-hyp} gives
\begin{equation}\label{eq:mu-size}
\mu
\ge
(4N)^{100(\gamma_*-\eps_*)}
\ge
(4N)^{10}.
\end{equation}
Hence $\mu\ge 2N$, and the avalanche hypotheses hold for every block family occurring at scales $N$ and $2N$.

Fix $X\in\{N,2N\}$ and $0\le r<n$. Set
\[
m_r=m_r(X):=\left\lfloor\frac{X-r}{n}\right\rfloor,
\qquad
s_r=s_r(X):=X-r-m_rn,
\qquad 0\le s_r<n,
\]
and
\[
R_r:=M_r(z),
\qquad
Q_{r,X}:=M_{m_rn}(T^rz),
\qquad
L_{r,X}:=M_{s_r}(T^{r+m_rn}z).
\]
Then
\begin{equation}\label{eq:factorization}
M_X(z)=L_{r,X}Q_{r,X}R_r.
\end{equation}
Since $\det M_j=1$, the two singular values of $M_j$ are reciprocal. Hence $\norm{M_j^{-1}}=\norm{M_j}.$ Using $0\le h_j\le1$,
\[
\log\norm{L_{r,X}^{\pm1}}\le s_rS,
\qquad
\log\norm{R_r^{\pm1}}\le rS.
\]
From \eqref{eq:factorization},
\begin{equation}\label{eq:Q-vs-X}
\left|
\log\norm{Q_{r,X}}-\log\norm{M_X(z)}
\right|
\le(r+s_r)S
<2nS.
\end{equation}

Write
\[
U_j(k):=\log\norm{M_j(T^kz)}.
\]
The avalanche principle applied to the $m_r$ matrices
\[
M_n(T^rz),\ M_n(T^{r+n}z),\ldots,
M_n(T^{r+(m_r-1)n}z)
\]
gives
\[
\left|
\log\norm{Q_{r,X}}
+\sum_{j=1}^{m_r-2}U_n(r+jn)
-\sum_{j=0}^{m_r-2}U_{2n}(r+jn)
\right|
\le C\frac{m_r}{\mu}.
\]
Replace $\log\norm{Q_{r,X}}$ by $\log\norm{M_X(z)}$ using \eqref{eq:Q-vs-X}, sum over $0\le r<n$, and divide by $nXS$. Lemma~\ref{lem:offset-counting}, with $X$ in place of $N$, yields
\begin{equation}
\left|
 h_X(z)
 +\frac1X\sum_{k=n}^{X-2n}h_n(T^kz)
 -\frac2X\sum_{k=0}^{X-2n}h_{2n}(T^kz)
\right|\le
C\left(
\frac nX
+
\frac{1}{nXS\mu}\sum_{r=0}^{n-1}m_r
\right).
\label{eq:truncated-averages}
\end{equation}
Since $m_r\le X/n$, we have
\[
\frac{1}{nXS\mu}\sum_{r=0}^{n-1}m_r
\le
\frac1{nS\mu}
\le
C\frac nX,
\]
where \eqref{eq:mu-size}, $X\le2N$, $n\ge1$, and $S\ge1$ were used in the last step.

The first average in \eqref{eq:truncated-averages} misses at most $3n$ indices from $\{0,\ldots,X-1\}$, the second misses at most $2n$. Since $0\le h_n,h_{2n}\le1$,
\begin{equation}\label{eq:fill-boundary}
\frac{3n}{X}+2\frac{2n}{X}
\le7\frac nX.
\end{equation}
Equations \eqref{eq:truncated-averages}--\eqref{eq:fill-boundary} prove \eqref{eq:APavg} simultaneously for $X=N$ and $X=2N$ on the same set $\mathcal B$.

For $X\in\{N,2N\}$, define
\[
\Phi_X(z)
:=
 h_X(z)
 +\frac1X\sum_{k=0}^{X-1}h_n(T^kz)
 -\frac2X\sum_{k=0}^{X-1}h_{2n}(T^kz).
\]
By \eqref{eq:APavg},
\[
|\Phi_X|\le C\frac nX
\quad\text{on }\mathcal B^c,
\qquad
|\Phi_X|\le4
\quad\text{on }\T^{2d}.
\]
Haar invariance gives
\[
\int_{\T^{2d}}\Phi_X
=
\ell_X+\ell_n-2\ell_{2n}.
\]
Hence, by \eqref{eq:Bbound},
\begin{equation}\label{eq:X-mean-bound}
|\ell_X+\ell_n-2\ell_{2n}|
\le
C\frac nX+C\mes\mathcal B
\le
C\left(\frac nN+Nq_n\right),
\qquad X\in\{N,2N\}.
\end{equation}
The case $X=N$ is \eqref{eq:APmean}. Subtracting the two inequalities in \eqref{eq:X-mean-bound} gives
\begin{equation}\label{eq:N-2N-abs}
|\ell_N-\ell_{2N}|
\le
C\left(\frac nN+Nq_n\right).
\end{equation}
Finally, subadditivity gives $L_{2N}\le L_N$, hence $\ell_{2N}\le\ell_N$. Combining with \eqref{eq:N-2N-abs} proves \eqref{eq:APgap}.
\end{proof}

\subsection{One-step propagation to the next scale}
\label{subsec:one-step-propagation}

The averaged avalanche comparison controls the large-scale quantity
$h_N$ up to orbit averages of $h_n$ and $h_{2n}$.  We now use the
short-block LDT from Proposition~\ref{prop:short-block-LDT} to replace
those orbit averages by their means outside a very small set.  This gives
a real-torus decomposition of $h_N-\ell_N$ into a uniformly small part
and an $L^1$ remainder.  Proposition~\ref{prop:boosting} then upgrades
that decomposition to a stretched-exponential LDT at scale $N$.

Fix the short-block constants $A,a,R_0$ from
Proposition~\ref{prop:short-block-LDT}. Increase $A$, if necessary, so that
$A\in\mathbb N$, $A\ge4$, and
\begin{equation}\label{eq:theta-def}
a<\theta:=1-\frac1A.
\end{equation}
Set
\begin{equation}\label{eq:fixed-parameters}
s:=\frac a2,
\qquad
P:=P_{2d}(s),
\qquad
b:=\sigma:=\beta_{2d}(s).
\end{equation}
Then $0<b=\sigma<1$.

\begin{proposition}\label{prop:step}
There are constants $C,c,N_1$, depending only on
$d,\rho,V,\gamma,\tau$ and the already fixed auxiliary exponents, with the
following property. They are independent of $n,N,\lambda,E$.
Suppose \eqref{eq:APscalebounds} holds and
\begin{equation}\label{eq:step-input}
N\ge \max\{N_1,R_0(2n)^A\},
\qquad
nS\ge100\log(4N),
\qquad
q_n\le(4N)^{-(P+4)}.
\end{equation}
Then
\begin{align}
\mes\bigl\{|h_N-\ell_N|>CN^{-b}\bigr\}
&\le Ce^{-cN^\sigma},\label{eq:step-LDT}\\
|\ell_N+\ell_n-2\ell_{2n}|
&\le C\frac nN,\label{eq:step-mean}\\
0\le \ell_N-\ell_{2N}
&\le C\frac nN.\label{eq:step-gap}
\end{align}
\end{proposition}

\begin{proof}
We first extract the mean and two-scale gap estimates from the averaged
avalanche comparison.  We then control the orbit averages of $h_n$ and
$h_{2n}$ by Proposition~\ref{prop:short-block-LDT}, and finally apply the
several-variable boosting estimate to $h_N-\ell_N$.

From $N\ge R_0(2n)^A$,
\begin{equation}\label{eq:step-scale-comparison}
n\le C_{A,R_0}N^{1/A},
\qquad
\frac nN\le C_{A,R_0}N^{-\theta}.
\end{equation}
After increasing $N_1$, $N\ge4n.$ Moreover,
\[
Nq_n\le 4^{-(P+4)}N^{-(P+3)}\le C N^{-1}
\le C\frac nN.
\]
Hence Lemma~\ref{lem:APaverage} gives
\begin{align}
|\ell_N+\ell_n-2\ell_{2n}|
&\le C\left(\frac nN+Nq_n\right)
\le C\frac nN,\label{eq:step-mean-proof}\\
0\le \ell_N-\ell_{2N}
&\le C\left(\frac nN+Nq_n\right)
\le C\frac nN.\label{eq:step-gap-proof}
\end{align}
This proves \eqref{eq:step-mean}--\eqref{eq:step-gap}.

Set
\[
A_n(z):=\frac1N\sum_{j=0}^{N-1}h_n(T^jz),
\qquad
A_{2n}(z):=\frac1N\sum_{j=0}^{N-1}h_{2n}(T^jz).
\]
By Lemma~\ref{lem:analytic-normalization}, $h_n$ and $h_{2n}$ satisfy the
hypotheses of Proposition~\ref{prop:short-block-LDT}. Since
$N\ge R_0(2n)^A\ge R_0n^A$,
\[
\mes\{|A_n-\ell_n|>CN^{-a}\}
+
\mes\{|A_{2n}-\ell_{2n}|>CN^{-a}\}
\le Ce^{-cN^a}.
\]

Let $\mathcal B$ be the bad set in Lemma~\ref{lem:APaverage}. Then $\mes\mathcal B\le CNq_n,$ and on $\mathcal B^c$,
\[
\left|h_N+A_n-2A_{2n}\right|
\le C\frac nN.
\]
Define
\begin{equation}\label{eq:step-G}
G:=\mathcal B^c
\cap\{|A_n-\ell_n|\le CN^{-a}\}
\cap\{|A_{2n}-\ell_{2n}|\le CN^{-a}\}.
\end{equation}
Then
\begin{equation}\label{eq:step-G-measure}
\mes(G^c)\le CNq_n+Ce^{-cN^a},
\end{equation}
and, for $z\in G$,
\[
|h_N(z)+\ell_n-2\ell_{2n}|
\le C\left(N^{-a}+\frac nN\right).
\]
Subtracting \eqref{eq:step-mean-proof},
\[
|h_N(z)-\ell_N|
\le C\left(N^{-a}+\frac nN\right).
\]
By \eqref{eq:theta-def} and \eqref{eq:step-scale-comparison}, after increasing
$N_1$,
\begin{equation}\label{eq:step-real-close}
\frac nN\le C N^{-\theta}\le C N^{-a},
\qquad
\sup_G|h_N-\ell_N|\le CN^{-a}.
\end{equation}

Put
\begin{equation}\label{eq:step-decomposition}
f:=h_N-\ell_N,
\qquad
f_0:=f\one_G,
\qquad
f_1:=f\one_{G^c}.
\end{equation}
Since $0\le h_N\le1$ and $0\le\ell_N\le1$ on the real torus,
$|f|\le1$. From \eqref{eq:step-real-close} and $s=a/2$,
\begin{equation}\label{eq:step-f0}
\norm{f_0}_{L^\infty}\le CN^{-a}\le N^{-s}
\end{equation}
for $N\ge N_1$. Furthermore, by \eqref{eq:step-G-measure} and
\eqref{eq:step-input},
\[
\norm{f_1}_{L^1}\le CNq_n+Ce^{-cN^a}\le C4^{-(P+4)}N^{-(P+3)}+Ce^{-cN^a}\le N^{-P},
\]
after one final increase of $N_1$.

By Lemma~\ref{lem:analytic-normalization}, $h_N$ is continuous and separately
subharmonic on
\[
\T^d_{\rho/(4N)}\times\T^d_{\rho/4},
\qquad
0\le h_N\le1.
\]
Hence $f=h_N-\ell_N$ is continuous and separately subharmonic there and
$|f|\le1$. Restricting the $y$-strip,
\begin{equation}\label{eq:step-isotropic-strip}
f:\T^{2d}_{\rho/(4N)}\longrightarrow\mathbb R,
\qquad |f|\le1.
\end{equation}
Apply Proposition~\ref{prop:boosting} with
\[
m=2d,
\qquad
\rho'=\frac\rho4,
\qquad
s=\frac a2,
\qquad
P=P_{2d}(s).
\]
Using \eqref{eq:step-f0}--\eqref{eq:step-isotropic-strip},
\begin{equation}\label{eq:step-boosted}
\mes\{|f|>CN^{-\beta_{2d}(s)}\}
\le Ce^{-cN^{\beta_{2d}(s)}}.
\end{equation}
Since $b=\sigma=\beta_{2d}(s)$, \eqref{eq:step-boosted} is
\eqref{eq:step-LDT}.
\end{proof}

% Section 7
\section{Proof of the Main theorem}
\label{sec:main-proof}

\subsection{All integer scales and the proof of positivity}
\label{sec:induction}
All constants in this section are fixed before the induction starts. The
parameters $A,a,R_0$ come from Proposition~\ref{prop:short-block-LDT}.
$P,b,\sigma,\theta$ are fixed by
\eqref{eq:theta-def}--\eqref{eq:fixed-parameters}.

Choose an integer $K\ge4$ such that
\begin{equation}\label{eq:K-choice}
R_0\left(\frac2K\right)^A\le1.
\end{equation}
Let $C_*$ dominate the constants in \eqref{eq:step-mean} and
\eqref{eq:step-gap}, and fix
\begin{equation}\label{eq:D-choice}
D\ge2C_*.
\end{equation}
Choose $\varepsilon_0>0$ so that
\begin{equation}\label{eq:eps0-choice}
2\varepsilon_0<\epss,
\qquad
6\varepsilon_0\le\frac1{16}.
\end{equation}
Choose
\begin{equation}\label{eq:Q-choice}
Q>A^2(P+4)+1.
\end{equation}
Apply Lemma~\ref{lem:initial-scale} with $(\varepsilon_0,Q)$ and denote by
$B>0$ and $\lambda_*>1$ the resulting constants.

For $n\ge2$, set
\begin{equation}\label{eq:base-scale-functions}
M(n):=\left\lceil(Kn)^A\right\rceil,
\qquad
B_0(n):=\left\lceil(2K M(n))^A\right\rceil.
\end{equation}
Then
\begin{equation}\label{eq:base-growth}
M(n)\le C_Kn^A,
\qquad
B_0(n)\le C_Kn^{A^2}.
\end{equation}
Hence \eqref{eq:Q-choice} implies, for all sufficiently large $n$,
\[
n^{-Q}\le \bigl(4B_0(n)\bigr)^{-(P+4)}.
\]

Choose $M_{\mathrm{ind}}$ so large that, for every $t\ge M_{\mathrm{ind}}$,
\begin{align}
Ct^{-b}&<\epss,
&
Dt^{-\theta}&<\frac{\gs}{40},\label{eq:terminal-smallness}\\
2Ce^{-ct^\sigma}
&\le \bigl(4(K(t+1))^A\bigr)^{-(P+4)},
&
t&\ge100\log\bigl(4(K(t+1))^A\bigr).
\label{eq:terminal-LDT-log}
\end{align}
Such a threshold exists. Choose $n_0$ so large that, with
\begin{equation}\label{eq:baseinterval}
M:=M(n_0),
\qquad
B_0:=B_0(n_0),
\end{equation}
one has
\begin{gather}
M\ge\max\{N_1,M_{\mathrm{ind}}\},
\qquad
n_0\ge100\log(4B_0),\label{eq:base-threshold-1}\\
n_0^{-Q}\le(4B_0)^{-(P+4)},
\qquad
C_*\frac{n_0}{M}\le\frac1{16},\label{eq:base-threshold-2}\\
\frac{3DM^{-\theta}}{1-2^{-\theta}}\le\frac18.
\label{eq:sumloss}
\end{gather}
All conditions are eventual as $n_0\to\infty$. Finally set
\begin{equation}\label{eq:lambda0-choice}
\lambda_0:=\max\{e^2,\lambda_*n_0^B\}.
\end{equation}
The constants $K,D,Q,B,n_0,M,B_0,\lambda_0$ are now fixed and are not
changed in the induction.

\begin{lemma}\label{lem:baseinterval}
Uniformly for $\lambda\ge\lambda_0$ and $E\in\mathbb R$, every integer
$N\in[M,B_0]$ satisfies
\begin{align}
\mes\bigl\{|h_N-\ell_N|>CN^{-b}\bigr\}
&\le Ce^{-cN^\sigma},\label{eq:indLDT}\\
\ell_N&\ge\frac78,\label{eq:basepositive}\\
0\le \ell_N-\ell_{2N}
&\le DN^{-\theta}.\label{eq:indgap}
\end{align}
\end{lemma}

\begin{proof}
Lemma~\ref{lem:initial-scale} and \eqref{eq:lambda0-choice} give
\begin{align}
&\mes\{|h_{n_0}-\ell_{n_0}|>2\varepsilon_0\}
+
\mes\{|h_{2n_0}-\ell_{2n_0}|>2\varepsilon_0\}
\le n_0^{-Q},\label{eq:base-initial-LDT}\\
&\ell_{n_0},\ell_{2n_0}\ge1-2\varepsilon_0,
\qquad
0\le\ell_{n_0}-\ell_{2n_0}\le2\varepsilon_0.
\label{eq:base-initial-means}
\end{align}
By \eqref{eq:eps0-choice},
\begin{equation}\label{eq:base-q}
q_{n_0}\le n_0^{-Q},
\end{equation}
and
\begin{equation}\label{eq:base-AP-hyp}
\ell_{n_0},\ell_{2n_0}>\gs,
\qquad
0\le\ell_{n_0}-\ell_{2n_0}<\frac{\gs}{40}.
\end{equation}

Fix $N\in[M,B_0]$. By \eqref{eq:K-choice} and
\eqref{eq:baseinterval},
\[
N\ge M\ge(Kn_0)^A\ge R_0(2n_0)^A.
\]
By $S\ge1$ and \eqref{eq:base-threshold-1},
\[
n_0S\ge n_0\ge100\log(4B_0)\ge100\log(4N).
\]
Finally, \eqref{eq:base-q}, \eqref{eq:base-threshold-2}, and $N\le B_0$
give
\[
q_{n_0}
\le n_0^{-Q}
\le(4B_0)^{-(P+4)}
\le(4N)^{-(P+4)}.
\]
Thus Proposition~\ref{prop:step} applies with short scale $n_0$ and long
scale $N$. It yields \eqref{eq:indLDT} and
\begin{align}
|\ell_N+\ell_{n_0}-2\ell_{2n_0}|
&\le C_*\frac{n_0}{N},\label{eq:base-mean-relation}\\
0\le\ell_N-\ell_{2N}
&\le C_*\frac{n_0}{N}.\label{eq:base-gap-raw}
\end{align}
Using \eqref{eq:base-initial-means} and \eqref{eq:base-threshold-2},
\begin{align}
\ell_N
&\ge 2\ell_{2n_0}-\ell_{n_0}-C_*\frac{n_0}{N}\notag\\
&=\ell_{n_0}-2(\ell_{n_0}-\ell_{2n_0})-C_*\frac{n_0}{N}\notag\\
&\ge1-6\varepsilon_0-C_*\frac{n_0}{M}\geq\frac78
\label{eq:base-positive-proof}
\end{align}
Moreover, $N\ge M\ge(Kn_0)^A$ implies
\[
n_0\le K^{-1}N^{1/A},
\qquad
\frac{n_0}{N}\le K^{-1}N^{-\theta}.
\]
Hence, by \eqref{eq:D-choice} and \eqref{eq:base-gap-raw},
\begin{equation}\label{eq:base-gap-proof}
0\le\ell_N-\ell_{2N}
\le \frac{C_*}{K}N^{-\theta}
\le DN^{-\theta}.
\end{equation}
The estimate at $2N$ in \eqref{eq:base-gap-raw} is part of the one-scale
conclusion at scale $N$, no assumption that $2N\in[M,B_0]$ is used.
\end{proof}

\begin{proposition}\label{prop:strong-induction}
Uniformly for $\lambda\ge\lambda_0$ and $E\in\mathbb R$, every integer
$N\ge M$ satisfies the LDT \eqref{eq:indLDT}, the gap estimate
\eqref{eq:indgap}, and
\begin{equation}\label{eq:allscale-positive}
\ell_N\ge\frac34.
\end{equation}
\end{proposition}

\begin{proof}
The interval $M\le N\le B_0$ is Lemma~\ref{lem:baseinterval}. Assume
$N>B_0$ and that the conclusion holds for every integer $j$ with
$M\le j<N$.

Set
\begin{equation}\label{eq:parent}
n=p(N):=\floor{\frac{N^{1/A}}{K}}.
\end{equation}
Since $N>B_0\ge(2KM)^A$,
$N^{1/A}/K>2M\ge2$. Hence
\begin{equation}\label{eq:parentbounds}
M\le n,
\qquad
\frac1{2K}N^{1/A}\le n\le\frac1K N^{1/A}.
\end{equation}
Moreover,
\begin{equation}\label{eq:2nlessN}
2n\le\frac2K N^{1/A}\le\frac12N^{1/A}\le\frac N2<N,
\end{equation}
and
\[
R_0(2n)^A
\le R_0\left(\frac2K\right)^AN
\le N.
\]
Finally,
\begin{equation}\label{eq:Nupper}
N<\bigl(K(n+1)\bigr)^A.
\end{equation}

By \eqref{eq:parentbounds}--\eqref{eq:2nlessN},
$M\le n<2n<N$. The induction hypothesis and
\eqref{eq:terminal-smallness} give
\begin{equation}\label{eq:inheritedgap}
\ell_n,\ell_{2n}\ge\frac34>\gs,
\qquad
0\le\ell_n-\ell_{2n}
\le Dn^{-\theta}\le\frac{\gs}{40}.
\end{equation}
Since $Cn^{-b}<\epss$, the inherited LDTs imply
\[
q_n:=\mes\{|h_n-\ell_n|>\epss\}
+\mes\{|h_{2n}-\ell_{2n}|>\epss\}
\le2Ce^{-cn^\sigma}.
\]
By \eqref{eq:terminal-LDT-log} and \eqref{eq:Nupper},
\[
q_n\le(4N)^{-(P+4)}.
\]
Also,
\[
nS\ge n
\ge100\log\bigl(4(K(n+1))^A\bigr)
\ge100\log(4N).
\]
Thus Proposition~\ref{prop:step} applies at the pair $(n,N)$. It yields
\eqref{eq:indLDT} at scale $N$ and
\begin{align}
|\ell_N+\ell_n-2\ell_{2n}|
&\le C_*\frac nN,\label{eq:meanrelation}\\
0\le\ell_N-\ell_{2N}
&\le C_*\frac nN.\label{eq:newgapraw}
\end{align}
By \eqref{eq:parentbounds},
\[
\frac nN\le\frac1K N^{-\theta},
\]
so
\begin{equation}\label{eq:newgap}
0\le\ell_N-\ell_{2N}
\le\frac{C_*}{K}N^{-\theta}
\le DN^{-\theta}.
\end{equation}

From \eqref{eq:meanrelation}, \eqref{eq:inheritedgap}, and
\eqref{eq:newgap},
\begin{align}
\ell_N
&\ge2\ell_{2n}-\ell_n-C_*\frac nN\notag\\
&=\ell_n-2(\ell_n-\ell_{2n})-C_*\frac nN\notag\\
&\ge\ell_n-2Dn^{-\theta}-DN^{-\theta}.
\label{eq:onesteploss}
\end{align}
Starting from $N_0=N$, define $N_{j+1}=p(N_j)$ while $N_j>B_0.$ The chain terminates at some $N_J\in[M,B_0]$ and $N_j\ge2N_{j+1}$ for $0\le j<J$. Applying \eqref{eq:onesteploss} at $N_0,\ldots,N_{J-1}$,
\[
\ell_N\ge\ell_{N_J}
-\sum_{j=0}^{J-1}
\bigl(2DN_{j+1}^{-\theta}+DN_j^{-\theta}\bigr)\ge\ell_{N_J}-3D\sum_{j=0}^{J-1}N_{j+1}^{-\theta}.
\]
Reading the chain upward from $N_J$,
\[
N_{J-r}\ge2^rN_J\ge2^rM,
\qquad0\le r\le J.
\]
Hence, by Lemma~\ref{lem:baseinterval} and \eqref{eq:sumloss},
\[
\ell_N\ge\frac78-3D\sum_{r=0}^{\infty}(2^rM)^{-\theta}=\frac78-\frac{3DM^{-\theta}}{1-2^{-\theta}}\ge\frac34.
\]
This closes the strong induction.
\end{proof}

\begin{proof}[Proof of \eqref{eq:mainLDT} and \eqref{eq:mainLE}]
Proposition~\ref{prop:strong-induction}, multiplied by $S$, proves
\eqref{eq:mainLDT} for every $N\ge M$. Moreover,
\begin{equation}\label{eq:LE-lower-final}
L(\lambda,E)=\lim_{N\to\infty}L_N(\lambda,E)
\ge\frac34S(\lambda,E)
\ge\frac12\log\lambda.
\end{equation}
Thus \eqref{eq:mainLE} follows from the all-scale lower bound, not merely
from positivity at the initial scale. Finally, dyadic summation of
\eqref{eq:indgap} gives
\begin{equation}\label{eq:convergence}
0\le L_N-L
=\sum_{j=0}^{\infty}(L_{2^jN}-L_{2^{j+1}N})
\le\frac{DS}{1-2^{-\theta}}N^{-\theta}.
\end{equation}
\end{proof}

\subsection{Continuity in the energy}
\label{sec:continuity}
The dyadic estimate \eqref{eq:convergence} gives
$L_N-L=O(SN^{-\theta})$. Combining it with the all-offset relation at an
exponentially longer scale yields the two-scale approximation needed for
energy continuity.

\begin{lemma}\label{lem:limitapprox}
There exist $C,c>0$ and $n_8\ge1$, independent of $E\in\mathbb R$ and
$\lambda\ge\lambda_0$, such that for every $n\ge n_8$,
\begin{equation}\label{eq:limitapprox}
|L(\lambda,E)+L_n(\lambda,E)-2L_{2n}(\lambda,E)|
\le CS(\lambda,E)e^{-cn^\sigma}.
\end{equation}
\end{lemma}

\begin{proof}
Write $h_m=u_m/S$ and $\ell_m=L_m/S$. By
Proposition~\ref{prop:strong-induction},
\[
\mes\{|h_m-\ell_m|>Cm^{-b}\}\le Ce^{-cm^\sigma},
\qquad
\ell_m\ge\frac34,
\qquad
0\le\ell_m-\ell_{2m}\le Dm^{-\theta}.
\]
After increasing $n_8$,
\[
Cn^{-b}<\epss,
\qquad
Dn^{-\theta}<\frac{\gs}{40},
\]
and hence
\begin{equation}\label{eq:limit-qn}
q_n:=\mes\{|h_n-\ell_n|>\epss\}
+\mes\{|h_{2n}-\ell_{2n}|>\epss\}
\le Ce^{-c_0n^\sigma}.
\end{equation}
Fix
\begin{equation}\label{eq:long-scale-choice}
0<\alpha<\frac{c_0}{4},
\qquad
N:=\left\lfloor e^{\alpha n^\sigma}\right\rfloor.
\end{equation}
For all sufficiently large $n$,
\[
N\ge4n,
\qquad
nS\ge n\ge100\log(4N),
\qquad
\ell_n,\ell_{2n}\ge\frac34>\gs,
\]
and the gap bound required by Lemma~\ref{lem:APaverage} follows from the
choice of $n_8$. Therefore
\begin{equation}
|L_N+L_n-2L_{2n}|\le CS\left(\frac nN+Nq_n\right)\le CSe^{-c_1n^\sigma},\label{eq:finite-two-scale}
\end{equation}
where
\[
\frac nN\le2ne^{-\alpha n^\sigma}
\le Ce^{-\alpha n^\sigma/2},
\qquad
Nq_n\le Ce^{-(c_0-\alpha)n^\sigma}.
\]
By \eqref{eq:convergence} and \eqref{eq:long-scale-choice},
\begin{equation}\label{eq:long-to-limit}
0\le L_N-L
\le\frac{DS}{1-2^{-\theta}}N^{-\theta}
\le CSe^{-\alpha\theta n^\sigma/2}.
\end{equation}
Combining \eqref{eq:finite-two-scale} and \eqref{eq:long-to-limit} proves
\eqref{eq:limitapprox}.
\end{proof}

\begin{proof}[Proof of \eqref{eq:continuitymain}]
Fix $\lambda\ge\lambda_0$ and a compact interval $I\subset\mathbb R$. Set
\begin{equation}\label{eq:BI-SI}
B_I:=\max\left\{2,\sup_{z\in\T^{2d},\,E\in I}
\norm{A(z;\lambda,E)}\right\},
\qquad
S_I:=\sup_{E\in I}S(\lambda,E)<\infty.
\end{equation}
For $E,E'\in I$, put $h:=|E-E'|$. Since
$\norm{A(z;\lambda,E)-A(z;\lambda,E')}=h$, telescoping the product gives
\[
\norm{M_n(z;E)-M_n(z;E')}
\le nB_I^{n-1}h.
\]
Moreover, $\det M_n=1$, hence
$\norm{M_n(z;E)},\norm{M_n(z;E')}\ge1$. Therefore
\begin{align}
|u_n(z;E)-u_n(z;E')|
&\le\frac1n
\left|\norm{M_n(z;E)}-\norm{M_n(z;E')}\right|\notag\\
&\le\frac1n\norm{M_n(z;E)-M_n(z;E')}\notag\\
&\le B_I^{n-1}h.
\end{align}
After integration,
\begin{equation}\label{eq:energyfinite}
|L_n(E)-L_n(E')|\le B_I^{n-1}|E-E'|.
\end{equation}
Applying Lemma~\ref{lem:limitapprox} at $E$ and $E'$, and then
\eqref{eq:energyfinite} at $n$ and $2n$,
\begin{equation}\label{eq:energy-balance}
|L(E)-L(E')|
\le C_Ie^{-cn^\sigma}+C_IB_I^{2n}h.
\end{equation}
For $h$ sufficiently small, choose
\[
n:=\left\lfloor\frac{|\log h|}{4\log B_I}\right\rfloor.
\]
After decreasing a fixed $h_0=h_0(I,\lambda)>0$,
\[
n\ge n_8,
\qquad
n\ge\frac{|\log h|}{8\log B_I},
\qquad
0<h<h_0.
\]
Hence
\[
B_I^{2n}h\le h^{1/2},
\qquad
e^{-cn^\sigma}\le\exp\bigl(-c_I|\log h|^\sigma\bigr).
\]
Since $0<\sigma<1$,
\[
h^{1/2}=e^{-|\log h|/2}
\le\exp\bigl(-c_I|\log h|^\sigma\bigr)
\]
for $h_0$ sufficiently small. Thus
\[
|L(\lambda,E)-L(\lambda,E')|
\le C_I\exp\bigl(-c_I|\log|E-E'||^\sigma\bigr),
\qquad0<|E-E'|<h_0.
\]
Finally, $0\le L(\lambda,E)\le S_I$ on $I$. Enlarging $C_I$ extends the
same estimate to $0<|E-E'|<1/2$, proving \eqref{eq:continuitymain}. Thus Theorem \ref{thm:main} holds.
\end{proof}
\begin{remark}
The proof above separates the argument into a finite-scale analytic input
and a deterministic propagation mechanism.  Once the large-coupling
initial scale and the short-block large-deviation estimate are available,
the passage to all sufficiently large integer scales uses only the
all-offset avalanche comparison and the strong induction developed in
Sections~\ref{sec:AP} and \ref{sec:main-proof}.  No additional
arithmetic input is required at the later scales.

The resulting modulus of continuity
\[
\exp\!\bigl(-c|\log h|^\sigma\bigr),
\qquad 0<\sigma<1,
\]
is weaker than a H\"older modulus but stronger than every negative power
of \(|\log h|\).  We do not attempt here to optimize this modulus.
\end{remark}

\vskip5mm

\section*{Acknowledgments}
\addcontentsline{toc}{section}{Acknowledgments}
This work was supported by the NSFC (grant nos. 11571327 and 11971059).

\vskip2mm

\section*{Data availability}
\addcontentsline{toc}{section}{Data availability}
No data were used for the research described in this article.

\vskip2mm

\section*{Declaration of competing interest}
\addcontentsline{toc}{section}{Declaration of competing interest}
The authors declare that they have no known competing financial
interests or personal relationships that could have appeared to
influence the work reported in this paper.

\vskip2mm

\section*{Declaration of Generative AI Use}
\addcontentsline{toc}{section}{Declaration of Generative AI Use}

During the preparation of this work, the authors used ChatGPT
(OpenAI) to assist with language editing, organization, literature
cross-checking, and presentation of the manuscript. After using this
tool, the authors reviewed and edited the content as needed,
independently verified the mathematical arguments and cited references,
and take full responsibility for the content of the publication.

\end{document}